\documentclass[12pt]{amsart}
\usepackage{amssymb,amscd}
\usepackage{graphicx, color}

\newtheorem*{Whitney towers}{Theorem~\ref{Whitney towers}}
\newtheorem*{h-towers}{Theorems ~\ref{half} \& \ref{$(n)$-solvable}}
\newtheorem*{surgery curves}{Theorem~\ref{surgery curves}}
\newtheorem*{cg=0}{Theorem~\ref{vanish}}

\newtheorem{thm}{Theorem}[section]
\newtheorem{mthm}[thm]{Main Theorem}
\newtheorem{lem}[thm]{Lemma}

\newtheorem{cla}[thm]{Claim}
\theoremstyle{definition}
\newtheorem{defn}[thm]{Definition}

\newtheorem{prob}[thm]{Problem}
\newtheorem{note}[thm]{Note}

\numberwithin{equation}{section}
\numberwithin{figure}{section}

\newcommand{\x}{\times}

\newcommand{\Z}{\mathbb{Z}}

\newcommand{\C}{\mathbb{C}}
\newcommand{\Q}{\mathbb{Q}}

\newcommand{\R}{\mathbb{R}}

\def\yen{{\setbox0=\hbox{Y}Y\kern-.97\wd0\vbox{hrule height.lex width.98%
\wd0\kern.33ex\hrule height.lex width.98\wd0\kern.45ex}}}

\begin{document}

\title{Local moves on knots and products of knots II}
\author{Louis H. Kauffman and  Eiji Ogasa}      

\thanks{\hskip-4mm E-mail: kauffman@uic.edu\quad  
ogasa@mail1.meijigakuin.ac.jp \newline
}

\date{}


\begin{abstract} 
We use the terms, knot product and local move, as defined in the text of the paper. 
Let $n$ be an integer$\geqq3$. 
Let $\mathcal S_n$ be the set of simple spherical $n$-knots in $S^{n+2}$. 
Let $m$ be an integer$\geqq4$. 
We prove that the map $j:\mathcal S_{2m}\to\mathcal S_{2m+4}$ is bijective, 
where $j(K)=K\otimes$Hopf,  and Hopf denotes the Hopf link. 


Let $J$ and $K$ be 1-links in $S^3$. 
Suppose that $J$ is obtained from $K$ by a single pass-move, 
which is a local-move on 1-links.  
Let $k$ be a positive integer. 
Let $P\otimes^kQ$ denote the knot product $P\otimes\underbrace{Q\otimes...\otimes Q}_k$. 
We prove the following: 
The $(4k+1)$-dimensional submanifold $J\otimes^k{\rm Hopf}$ $\subset S^{4k+3}$
is obtained from 
$K\otimes^k{\rm Hopf}$  by a single $(2k+1,2k+1)$-pass-move,   
which is a local-move on $(4k+1)$-submanifolds contained in $S^{4k+3}$.  
See the body of the paper for the definitions of all local moves in this abstract. 

We prove the following: 
Let $a,b,a',b'$ and $k$ be positive integers. 
If the $(a,b)$ torus link is pass-move equivalent to the $(a',b')$ torus link, 
then the Brieskorn manifold 
$\Sigma(a,b,\underbrace{2,...,2}_{2k})$  
is diffeomorphic to $\Sigma(a',b',\underbrace{2,...,2}_{2k})$ 
as abstract manifolds.

Let $J$ and $K$ be (not necessarily connected or spherical) 2-dimensional closed oriented submanifolds  
in $S^4$. 
Suppose that $J$ is obtained from $K$ by a single ribbon-move, which is a local-move on 2-dimensional submanifolds contained in $S^4$.  
Let $k$ be an integer$\geq2$. 
We prove the following: 
The $(4k+2)$-submanifold  $J\otimes^k{\rm Hopf}$ $\subset S^{4k+4}$
is obtained from 
$K\otimes^k{\rm Hopf}$  by a single $(2k+1,2k+2)$-pass-move,   
which is a local-move on $(4k+2)$-dimensional submanifolds contained in $S^{4k+4}$.

\end{abstract} 
\maketitle

\section{Introduction} \label{Introduction}  
\noindent
The purpose of this paper is to investigate certain key constructions on a codimension-two differentiable closed oriented submanifold of a standard sphere $S^{n}$ of dimension $n$ in different dimensions. 
Of course  knots are examples of these submanifolds 
(see the next section for the definition of knots and that of submanifolds).      
\\

The first construction is the {\it knot product} 
(It is defined in \cite[Definnition 3.3]{KauffmanNeumann}.  
An alternative definition is given in (\cite[Lemma 3.4]{KauffmanNeumann}.)  
Let $k$ and $l$ be positive integers. 
For  
a (not necessarily connected or spherical) $(k-2)$-dimensional closed oriented submanifold $K$ of $S^k$
and 
a (not necessarily connected or spherical) $(l-2)$-dimensional 
closed oriented submanifold $L$ of $S^l$,    
a knot product  
$K\otimes L$ is defined 
as long as $K$ or $L$ is fibered. 
$K\otimes L$ is 
a (not necessarily connected or spherical) $(k+l-1)$-dimensional closed oriented submanifold of $S^{k+l+1}$. 
 (See precise definitions 
in the next section.  
$K\otimes L$ is defined even if $k$ or $l$ is one.) 
We can handle not only the case where $K$,  $L$, and  $K\otimes L$ are knots 
but also the case where not all of $K$,  $L$, and  $K\otimes L$  are knots 
although the name is `knot' products.

The knot product generalizes a property of the links of complex hypersurface singularities 
(see 
\cite[Proposition 4.3]{KauffmanNeumann}), 
and can be used to give constructions for Brieskorn (sub)manifolds 
(\cite[the sixth line from the bottom of page 1106]{Kauffman} and 
\cite[the last sentence of \S6, and the first remark in page 390]{KauffmanNeumann}),  
and other phenomena in low and high dimensional topology.   
\\

The second construction is a generalization of 
an equivalence relation in classical knot theory that we have called ``pass equivalence". 
Two classical knots (embeddings of the circle)
are said to be {\it pass equivalent} if there is a sequence of pass-moves as shown in the figure of Definition \ref{Komagome} combined with ambient isotopy, taking one knot to the other.
Two single component classical knots are pass equivalent if and only if they have the same Arf invariant 
(\cite[Proposition 5.6 in page 77 and Corollary in page 260]{Kauffmanon}).  
In this paper, we study high dimensional analogues of the pass move, 
which  \cite{Ogasa98n}  begun, 
and we show that high dimensional pass equivalence works compatibly in relation to taking the product with the Hopf link, which \cite{KauffmanOgasa} begun. 
Taking products with the Hopf link preserves Seifert linking pairings and is a natural operation in the context of knot products and high dimensional knot theory. 
It is important for us to know that the pass equivalence relation is preserved under this product.(The above property of Seifert pairings is proved in \cite[Proposition 6.2]{KauffmanNeumann}. 
\cite[Proposition 6.2]{KauffmanNeumann} is proved by using 
\cite[Lemma 6.1]{KauffmanNeumann}. 
We will explain  \cite[Lemma 6.1]{KauffmanNeumann} in more detail two paragraphs below.)
\\

In the case of high dimensional knot theory we still need to understand a natural framework for the many phenomena that occur there. The Brieskorn (sub)manifolds were the first construction that showed that exotic differentiable structures on spheres were natural in the sense that these 
spheres occurred as links of algebraic singularities. 
Not only did they occur as such links, but their exoticity was related to the way those spheres were knotted in codimension two (\cite[the last sentence of \S8]{Milnor}). 
The knot product construction  provides  
a systematic method 
for producing Brieskorn (sub)manifolds and their generalizations 
(\cite[the sixth line from the bottom of page 1106]{Kauffman} and 
\cite[the last sentence of \S6, and the first remark in page 390]{KauffmanNeumann}), 
and it provides a linkage between different dimensions of construction 
(see the previous paragraph, Main Theorems \ref{coffee}, \ref{sugar}, and \ref{melon}). 
Using knot products, one sees that Brieskorn manifolds are obtained by an iterated branched covering construction applied to torus knots. 
By making these new constructions for manifolds and knots in higher dimensions we open the possibilities of new understandings and new questions about the nature of higher dimensional space.
Remember that there remain many important open problems in high dimensional knot and link theory (\cite{CochranOrr, CochranOrrTeichner, Ogasanewpair}).  \\

Let $k$ and $l$ be positive integers. 
Let $K$ (respectively, $L$) be 
a codimension two closed oriented submanifold of $S^k$ (respectively $S^l$), 
and $J$ or $K$ fibered. 
A relation among a Seifert hypersurface for $J$, that for $K$, and that for $J\otimes K$ 
is shown in \cite[Lemma 6.1]{KauffmanNeumann}. 
The relation deduces a relation among 
a Seifert matrix for $J$, that for $K$, and that for $J\otimes K$ 
(\cite[Proposition 6.2]{KauffmanNeumann}).  
Indeed the relation which is shown in \cite[Lemma 6.1]{KauffmanNeumann}  
implies stronger facts. It is a crucial idea of this paper. 
See Claims \ref{dankai}, \ref{dansui}, \ref{Xdankai}, and \ref{Xdansui}, and 
Theorem \ref{mars}. These imply our main results
(Main Theorems \ref{coffee}, \ref{sugar}, and \ref{melon}). 
\\

In \cite{KauffmanOgasa} we began to research relations between knot products and 
local-moves on knots. 
This paper is a sequel to \cite{KauffmanOgasa}. 
\cite{KauffmanOgasaII} is a preprint of this paper. 
Knot products, or products of knots, were defined in \cite{KauffmanNeumann}
(It is defined in \cite[Definnition 3.3]{KauffmanNeumann}.  
An alternative definition is given in \cite[Lemma 3.4]{KauffmanNeumann}.)  
\cite{Kauffman} is an announcement of \cite{KauffmanNeumann}. 
Knot products have been researched in 
\cite{Kauffmansuper, KauffmanNeumann, KauffmanOgasa}
(See also \cite{revisited} as it uses the same ideas as \cite{Kauffmansuper} in dimensions 3 and 4.)  
Local moves on high dimensional knots were defined in \cite{Ogasa98n, Ogasa04,Ogasa09}
and 
have been  researched in 
\cite{KauffmanOgasa, Ogasa98n,  Ogasa02,  Ogasa04, Ogasa07,   Ogasa09, OgasaT3, OgasaIH, 
OgasaZAlex, Ogasaribbontwo}. 
\cite {KauffmanOgasaB} is a sequel of this paper. 

Our main results are Main Theorems \ref{coffee}, \ref{sugar}, and \ref{melon}.

\tableofcontents

\bigbreak
\section{A main result on knot products}\label{amain}
\noindent 
We work in the smooth category unless we indicate otherwise.

Let $n$ and $m$ be positive integers.  
An  {\it $($oriented$) ($ordered$)$ $m$-component $n$-$($dimensional$)$ 
link} 
 is a smooth, oriented submanifold $L=(L_1,...,L_m)$ $\subset$ $S^{n+2}$, 
which is the ordered disjoint union of $m$ manifolds, each PL homeomorphic 
to the standard $n$-sphere. If $m=1$, then $L$ is called a {\it 
knot}.  
See \cite[Chapter 1 section 1]{KobayashiNomizu} for the definition of smooth manifolds and that of smooth submanifolds. 


Note the following: 
We usually define $n$-links as above (see e.g. \cite{CochranOrr}). 
Not all $n$-knots are diffeomorphic to the standard $n$-sphere 
although all $n$-knots are PL homeomorphic to the standard $n$-sphere.  
Indeed many exotic $n$-spheres,  
which are not diffeomorphic to the standard $n$-sphere,  
can be embedded smoothly in $S^{n+2}$ (see \cite{Levinecob, Levinesimp, Milnor} for the proof).  
We exclude wild knots as usual (see \cite[1.3 Definition and the parts near it]{BZ}). 

In this paper, diffeomorphism (respectively, PL homeomorphism) means 
orientation preserving diffeomorphism (respectively, PL homeomorphism) unless otherwise stated. 
\\

Let $n$ be a positive integer. 
Two submanifolds $J$ and  $K$ $\subset S^n$ are {\it $($ambient$)$ isotopic}  
if there is a smooth orientation preserving family of diffeomorphisms $\eta_t$ of $S^n$, $0\leqq t\leqq1$, with $\eta_0$ the identity and $\eta_1(J)=K$.   
\\

An $m$-component $n$-link $L=(L_1,...,L_m)$ is called a {\it trivial $(n$-$)$link}   
if each $L_i$ bounds an $(n+1)$-ball $B_i$ embedded in $S^{n+2}$ 
and if $B_i\cap B_j=\phi$ for each distinct $i, j$. 
If $L$ is a trivial 1-component ($n$-)link,  then $L$ is called a {\it trivial $(n$-$)$knot}. 
\\

Let $x$ and $y$ be nonnegative integers.   
Let $X$ be an $x$-dimensional submanifold of a $y$-dimensional manifold $Y$. 
Let $N(X)$ be the tubular neighborhood of $X$ in $Y$. 
See \cite[page 115]{MilnorStasheff} for the definition and properties of tubular neighborhoods. 
In this paper the tubular neighborhood means 
the closed tubular neighborhood not the open tubular neighborhood.  
The fiber of the closed (respectively, open) tubular neighborhood is 
the closed (respectively, open) $(y-x)$-ball.  
We do not say what $Y$ is, if there is no danger of confusion.

Let $n$ be a positive integer. 
Let $K$ be 
a (not necessarily connected or spherical) $n$-dimensional closed oriented submanifold of $S^{n+2}$. 
Let $N(K)$ be the tubular neighborhood of $K$ in $S^{n+2}$. 
Note that $\partial N(K)$ is the total space of a product bundle $\partial D^2\x K\to \partial D^2$.   
Let $\overline{S^{n+2}-N(K)}$ be the closure of $S^{n+2}-N(K)$ in $S^{n+2}$. 
We say that $K$ is a {\it fibered submanifold} 
if $\overline{S^{n+2}-N(K)}$ is the total space of a fiber bundle $\pi:\overline{S^{n+2}-N(K)}\to S^1$  
and 
if $\pi|_{\partial N(K)}:\partial N(K)\to S^1$ is the projection map of a product fiber bundle 
$\partial D^2\x K\to \partial D^2$.  
Here, note that $\partial D^2\x K=\partial N(K)$ and $\partial D^2=S^1$.  
\\

Let $l$ be a positive integer.  
Let $L$ be a (not necessarily connected or spherical) closed oriented $l$-submanifold contained in $S^{l+2}$. 
Let $W$ be a connected compact oriented $(l+1)$-submanifold of $S^{l+2}$ 
such that $\partial W=A$. 
We call $W$ a {\it Seifert hypersurface} for $A$.   
In this paper we abbreviate `manifold-with-boundary'  (respectively, `submanifold-with-boundary')  to `manifold' (respectively, `submanifold') when the meaning is clear from the context.  
\\

We review the definition of knot products, or products of knots, after we explain a question 
which 
motivates the definition of knot products.  
This question was discussed and answered in \cite[\S 4]{KauffmanNeumann}. 

Let $a$ be a positive integer. 
Let $f: \C^a \longrightarrow\C$ 
be a (complex) polynomial mapping with an isolated singularity at the origin of $\C^a$
and satisfy $f$(the origin)=the origin. 
Recall that  one defines 
the {\it link of this singularity}, $L(f)=f^{-1}(0)\cap S^{2a-1}\subset S^{2a-1}.$ 
(See  \cite{Milnor} for the precise definition of the link of singularity.)
Here the symbol $f^{-1}(0)$ denotes the variety of $f$,  
and $S^{2a-1}$ denotes a sufficiently small sphere about the origin of $\C^a.$

Let $b$ be a positive integer. Given another such polynomial $g:\C^{b} \longrightarrow \C$, 
 form $f + g: \C^{a+b} =\C^a \times \C^b\to \C$ 
by $(f+g)(x,y)=f(x)+g(y)$. 

What kind of relation do we have among  
$L(f) \subset S^{2a-1}$, $L(g) \subset S^{2b-1}$, and $L(f + g)\\\subset S^{2a+2b-1}$?  
An answer to this question is explained 
below in the form of a construction that works for more general codimension two embeddings. 
This construction is due to \cite{KauffmanNeumann}. 
%
\\

Knot products are defined in \cite[Definition 3.3]{KauffmanNeumann}. 
We have described them in the second paragraph of \S\ref{Introduction} of this paper. 
An alternative definition is given in \cite[Lemma 3.4]{KauffmanNeumann}, 
which we cite almost verbatim below. 
Some notations in the cited lemma 
 are different from ones in this paper. See the remarks below the citation.

We emphasize that the reference  \cite{KauffmanNeumann} contains the proofs about 
the well-definedness of this construction.

A knot $(S^k, K)$ 
is defined in \cite[the fourth paragraph of page 371]{KauffmanNeumann} as   
a (not necessarily connected or spherical) $(k-2)$-dimensional closed oriented submanifold $K$ of $S^k$. 
See the footnote
\footnote{ 
As \cite{KauffmanNeumann, Milnor} etc. do, we sometimes use the term, `knots', for   
 (not necessarily connected or spherical) closed oriented submanifolds. 
Then we sometimes use the term, `spherical knots', 
for what we defined in the second paragraph of this section. 


Recall the following well-known fact: 
Let $n$ be a positive integer$\neq4$. 
Let $M$ be a smooth manifold. 
$M$ is homeomorphic to  the standard $n$-sphere if and only if 
$M$ is PL homeomorphic to  the standard $n$-sphere. 
If $n=4$, it is an open problem whether this fact holds. 
So, in the $n\neq4$ case,  
we should use both terms, 
 PL spherical knots and topological spherical knots as for now.    
Of course we do not use these terms when the meaning is clear from the context. 

Note that we always work in the differentiable category in the present paper.   
}
for the term, `knot'.

A fibered knot $(S^l, L, b)$ 
is defined in \cite[the fourth paragraph of page 371]{KauffmanNeumann} as   
a (not necessarily connected or spherical) $(l-2)$-dimensional fibered closed oriented submanifold $L$ of $S^l$. 

The key to the knot product construction is the association of maps to $D^2$ (as explained below) 
with a fibered closed oriented submanifold $(S^l,L,b)$. 
The fibered submanifold has a map $b:\overline{S^l-N(L)}\to S^1$
that is a fibration and restricts to the projection to $S^1$ on the boundary 
of $N(L)\cong L\x D^2$, where $\cong$ denotes an orientation preserving diffeomorphism. 
Thus $b$ extends to a map (still called) $b:S^l\to D^2$ with $b^{-1}(0)=L$. We can then 
form the cone on this map, $cb:D^{l+1}\to D^2$ with $(cb)^{-1}(0)=CL$, 
where $CL$ denotes the cone on $L$. 
Define a pair $(D^{2a}, f^{-1}(0))$ by using the algebraic variety several paragraphs above, 
where $D^{2a}$ is a $2a$-ball differentiably embedded in $\C^{2a}$ whose boundary is $S^{2a-1}$. 
The pair $(D^{l+1}, CL)$ is 
a generalization of $(D^{2a}, f^{-1}(0))$   
to the  case where we do not use an algebraic variety. 
\\

\noindent{\bf Definition of knot products.} 
Let $k$ and $l$ be positive integers.  Let $S^k=\partial D^{k+1}$ and $S^l=\partial D^{l+1}$.

\noindent{\bf Lemma 3.4 of \cite{KauffmanNeumann}.} 
{\it 
Let $(S^k, K)$ be a knot 
and 
$(S^l,L,b)$ a fibered knot. 
Let $F\subset D^{k+1}$ be a spanning manifold for $K$ 
as 
in \cite[Definition 3.1]{KauffmanNeumann}. 
Use \cite[Lemma 2.3]{KauffmanNeumann} to obtain $\gamma:D^{k+1}\to D^2$ with 
$\gamma^{-1}(0)=F,0$ a regular value of $\gamma$. 
Let $\tau:D^{l+1}\to D^2$ be a smoothing of $cb$.  
Use these maps to form the pullback 

$$
\begin{CD}
b(D^{k+1},F)&@>>>&D^{l+1}\\
@VVV& &@VV{\tau}V\\
D^{k+1}&@>{\gamma}>>&D^2
  \end{CD}
$$
\\


\noindent 
Thus 
$b(D^{k+1},F)\subset D^{k+1}\x D^{l+1}$. 
Then the knot product $(S^{k+l+1}, K\otimes b)$ 
is obtained by taking boundaries from this embedding. 
That is,   
$(S^{k+l+1}, K\otimes b)$ is isotopic to \\
$(\partial(D^{k+1}\x D^{l+1}), \partial(b(D^{k+1},F)))$.  }

$(S^{k+l+1}, K\otimes b)$  is well-defined in the differentiable category in terms of the embedding $K$ and the fiber structure of $L$. 
\\

For 
$(S^k, K)$ and 
$(S^l, L, b)$, 
$(S^{k+l+1}, K\otimes b)$ is indicated by $(S^{k+l+1}, K\otimes L)$ as written in  
\cite[the first paragraph of \S 3]{KauffmanNeumann}. 



See \cite[Fig. 1 in page 371 and the part above it]{KauffmanNeumann} for 
the word `smoothing' in the cited lemma. 
In particular, a smoothing of $cb:D^{l+1}\to D^2$ (as in \cite[Lemma 3.4]{KauffmanNeumann}) 
with $(cb)^{-1}(0)=CL$ is obtained by changing to a nearby smooth map 
$\widetilde{cb}:D^{l+1}\to D^2$ with $(\widetilde{cb})^{-1}(0)=G$, 
where $G\subset D^{l+1}$ is a smooth submanifold of $D^{l+1}$ 
and $(\widetilde{cb})^{-1}(0)\cap S^{l}$ is isotopic to $L$. 
The manifold $G$ is part of a family of submanifolds 
that deforms to the singular subspace $CL$.

In this paper we let $K\otimes L$ denote 
$(S^{k+l+1}, K\otimes L)$ in order to make the notation shorter.  
Note that in 
\cite[the first paragraph of \S 3, and Definition 3.1]{KauffmanNeumann},  
$K\otimes L$ denotes the diffeomorphism type of the submanifold $(S^{k+l+1}, K\otimes L)$, 
$K\otimes L\hookrightarrow S^{k+l+1}$.   

If $k$ (respectively, $l$) is one, we regard $K$ (respectively $L$) as the empty knot, 
which is defined in \cite[line($-3$) of page 371]{KauffmanNeumann}. 
The empty knots correspond to the maps \\$\tau_a:D^2\to D^2, \tau_a(z)=z^a$, 
where $a$ is an integer$\geqq2$. 
Let $p$ be a nonnegative integer. 
Let $E$ be a $p$-dimensional $($not necessarily connected or spherical$)$ 
closed oriented submanifold contained in $S^{p+2}$. 
In Lemma \ref{soda} we explain 
how we construct 
the $(p+2)$-dimensional closed oriented submanifold 
$E\otimes[2]$ in the standard $(p+4)$-sphere:  
$E\otimes[a]$ is an embedding of 
the $a$-fold branched cyclic cover of $S^{p+2}$ along $E\subset S^{p+2}$ in $S^{p+4}$.    
See Lemma \ref{soda} for the detail.

$b(D^{k+1},F)$ is a submanifold of $D^{k+1}\x D^{l+1}$ 
as \cite[the first line of page 374]{KauffmanNeumann} is pointed out:  
%
In \cite[the last line and the diagram on it of page 373]{KauffmanNeumann},   
let $M$ be $D^{k+1}$, $D^{n+1}$ be $D^{l+1}$ and $V$ be $F$  
(note this $F$ is not `$F$ in \cite[Definition 2.1]{KauffmanNeumann}').   
Then $N=D^{k+1}\x D^{l+1}$. 
%
Note  
$N=\tau(M, \alpha)$ in \cite[the third line of Theorem 2.2]{KauffmanNeumann}, 
$\tau(M, \alpha)=b(M, \alpha)$ in \cite[Theorem 2.2.(ii)]{KauffmanNeumann}, and    
$b(M,\alpha)=b(M,V)$ in  \cite[Remark in page 376]{KauffmanNeumann}.  
Since $b(D^{k+1},F)$ is an inverse image of a regular value by a smooth map,  
$b(D^{k+1},F)$ is a submanifold of $D^{k+1}\x D^{l+1}$.

Hence $\partial(b(D^{k+1},F))$ 
is a well-defined differentiable submanifold of $S^{k+l+1}$.


The first remark in \cite[page 380]{KauffmanNeumann} explains the reason for the following: 
When we consider the knot product of two submanifolds,  
we impose the condition that one of the two submanifolds is fibered. 
If neither submanifold is fibered, the construction is not uniquely defined.

In Lemmas \ref{soda} and \ref{sodasui}, 
we show an example of \cite[Lemma 3.4]{KauffmanNeumann}. 
\\

\cite[Corollary 3.6]{KauffmanNeumann} implies that for fibered knots $(S^l,L,b)$ and $(S^{l'},L',b')$,  \\
$(S^{l+l'+1}, L\otimes L')$ is isotopic to $(-1)^{(l-1)(l'-1)}(S^{l+l'+1}, L'\otimes L)$. 
%
\\

\cite[Proposition 4.3]{KauffmanNeumann} gives an answer to the question 
which we posed several paragraphs above.  
It shows that $L(f)\otimes L(g)$ is isotopic to $L(f + g)$.   

Note that knot products $J\otimes K$ are defined 
even when neither $J$ nor $K$ is a link of singularity. 
This fact is significant because it utilizes structure that formerly was only available through algebraic varieties to knot theory proper. 


\bigbreak 
Let $n$ be a positive integer. 
Let $\mathcal E_n$ be 
the set of (not necessarily connected or spherical) $n$-dimensional closed oriented 
submanifolds in $S^{n+2}$. 
We abbreviate $A\otimes\text{(the Hopf link)}$ to $A\otimes$Hopf.  
Let $K\in\mathcal E_n$. 
By the definition of the knot product, 
which we review before here in this section, 
$K\otimes\mathrm{Hopf}\in\mathcal E_{n+4}$. 
Thus we obtain a map 
$$\otimes\mathrm{Hopf}: \mathcal E_{n}\to\mathcal E_{n+4}.$$
$$\hskip23mm K\mapsto K\otimes\mathrm{Hopf}.$$ 

Let $\mathcal A$ be a subset of $\mathcal E_n$.  
We abbreviate a map $(\otimes\mathrm{Hopf})\vert_{\mathcal A}$ 
to  $\otimes\mathrm{Hopf}$. 
Let 
$\otimes\mathrm{Hopf}(\mathcal A)$ 
be denoted by 
$\mathcal A\otimes\mathrm{Hopf}$. 



Let $p$ be a positive integer.       
Let $K$ be a spherical $(2p+1)$-knot contained in $S^{2p+3}$.  
Let $i$ be an integer such that $2\leqq i\leqq p$. 
We say that $K$ is a {\it simple spherical $(2p+1)$-knot} if $K$ satisfies that 
$\pi_1(S^{2p+3}-K)\cong\Z$ and 
$\pi_i(S^{2p+3}-K)\cong0$. 

Note that for each integer $j$, 
$
\pi_j(S^{2p+3}-{\rm Int}N(K))
\cong
\pi_j(S^{2p+3}-N(K))
\cong
\pi_j(S^{2p+3}-K). 
$ 
In knot theory we often use the compact manifold $S^{2p+3}-{\rm Int}N(K)$. 

See \cite{Farber1978, Farber1980, Farber1983, Farber1984I, Farber1984II, Levinesimp} 
for  important results on simple spherical knots. 
Let $n$ be any integer$\geqq3$. 
Let $\mathcal S_{n}$ be the set of simple spherical $n$-knots in $S^{n+2}$. 
The following were proved.

\begin{thm}\label{K} 
{\rm{\bf{(\cite[Theorem 6]{Kauffman}.)}}} 
Let $n$ be a positive integer.  
Let $\mathcal K_{n}$ be the set of $n$-dimensional spherical knots in $S^{n+2}$. 
Then we have the following.

\smallbreak
\noindent$(1)$ 
$\mathcal K_{n}\otimes\mathrm{Hopf}\subset\mathcal K_{n+4}.$

\smallbreak
\noindent$(2)$  
Let $K, J\in\mathcal K_{n}$. 
Suppose that $K$ is knot-cobordant to $J$. 
Then $K\otimes\mathrm{Hopf}$ is 

\noindent{\color{white}$(2)$}   
knot-cobordant to $J\otimes\mathrm{Hopf}$.  
Thus we can define a homomorphism  
$$ \otimes\mathrm{Hopf}: C_n\to C_{n+4},$$

\noindent{\color{white}$(2)$}   
where $C_n$ is the $n$-dimensional knot cobordism group. 

\smallbreak
\noindent$(3)$ 
$ \otimes\mathrm{Hopf}: C_n\to C_{n+4}$ is an isomorphism if $n\neq1,3$. 

\noindent{\color{white}$(3)$}   
$\otimes\mathrm{Hopf}: C_3\to C_7$ is injective and not surjective. 

\noindent{\color{white}$(3)$}   
$\otimes\mathrm{Hopf}: C_1\to C_5$ is surjective and not injective.

\smallbreak
\noindent$(4)$ 
$\mathcal S_{n}\otimes\mathrm{Hopf}\subset\mathcal S_{n+4}.$

\end{thm}

\smallbreak
\noindent
{\bf Note.}
See \cite{Kervaire, Levinecob} for the definition of knot-cobordism and 
that of the $n$-dimensional knot cobordism group $C_n$. 
\smallbreak

\begin{thm}\label{LE}   {\rm{\bf{(\cite[Theorem 11.7]{KauffmanOgasa}.)}}}  
We have the following. \newline
For  any integer $m\geqq2$, 
$\otimes\mathrm{Hopf}: \mathcal S_{2m+1}\to \mathcal S_{2m+5}$ 
is a bijective map.  
\newline
$ \otimes\mathrm{Hopf}: \mathcal S_3\to \mathcal S_7$ 
is injective and not surjective. \newline
$ \otimes\mathrm{Hopf}: \mathcal K_1\to \mathcal S_5$ is surjective and not injective. 
\end{thm}

In this paper we prove the following.

\begin{mthm} \label{coffee}    
For each integer $m\geqq4$, 
$\otimes\mathrm{Hopf}: \mathcal S_{2m}\to\mathcal S_{2m+4}$ is a bijective map.  
\end{mthm}

We will review some more notations in \S\ref{review} 
before we will state other main theorems in \S\ref{main2}.

\bigbreak
\section{Review of  the definitions and some results on local-moves on $n$-dimensional knots,  
where $n$ is a positive integer} \label{review}

\begin{defn}\label{Komagome} {\bf(\cite[Definition 5.4]{Kauffmanon}.)}
Two 1-links are {\it pass-move-equivalent}  if one is obtained from the other 
by a sequence of pass-moves. 
See 
the following figure 
 for an illustration of the pass-move. 
Each of four arcs in the 3-ball 
may belong to different components of the 1-link.

\vskip3mm
\hskip18mm
\unitlength 0.1in
\begin{picture}(36.39,15.20)(4.01,-22.91)
%
\special{pn 8}%
\special{ar 1138 1552 737 739  0.0000000 6.2831853}%
%
\special{pn 8}%
\special{pa 566 2017}%
\special{pa 1499 921}%
\special{fp}%
\special{sh 1}%
\special{pa 1499 921}%
\special{pa 1441 959}%
\special{pa 1464 962}%
\special{pa 1471 985}%
\special{pa 1499 921}%
\special{fp}%
%
\special{pn 8}%
\special{pa 1705 1085}%
\special{pa 771 2187}%
\special{fp}%
\special{sh 1}%
\special{pa 771 2187}%
\special{pa 829 2149}%
\special{pa 805 2146}%
\special{pa 799 2123}%
\special{pa 771 2187}%
\special{fp}%
%
\special{pn 8}%
\special{pa 1292 1510}%
\special{pa 1143 1401}%
\special{fp}%
%
\special{pn 8}%
\special{pa 1132 1707}%
\special{pa 988 1582}%
\special{fp}%
%
\special{pn 8}%
\special{pa 2375 1660}%
\special{pa 2091 1661}%
\special{fp}%
\special{sh 1}%
\special{pa 2091 1661}%
\special{pa 2158 1681}%
\special{pa 2144 1661}%
\special{pa 2158 1641}%
\special{pa 2091 1661}%
\special{fp}%
%
\special{pn 8}%
\special{pa 2102 1541}%
\special{pa 2401 1541}%
\special{fp}%
\special{sh 1}%
\special{pa 2401 1541}%
\special{pa 2334 1521}%
\special{pa 2348 1541}%
\special{pa 2334 1561}%
\special{pa 2401 1541}%
\special{fp}%
\put(18.7000,-21.7000){\makebox(0,0)[lb]{{\bf pass-move}}}%
%
\special{pn 8}%
\special{ar 3303 1510 737 739  0.0000000 6.2831853}%
%
\special{pn 8}%
\special{pa 3148 1504}%
\special{pa 3272 1360}%
\special{fp}%
%
\special{pn 8}%
\special{pa 3298 1733}%
\special{pa 2968 2167}%
\special{fp}%
\special{sh 1}%
\special{pa 2968 2167}%
\special{pa 3024 2126}%
\special{pa 3000 2125}%
\special{pa 2992 2102}%
\special{pa 2968 2167}%
\special{fp}%
%
\special{pn 8}%
\special{pa 3102 1582}%
\special{pa 2751 1991}%
\special{fp}%
%
\special{pn 8}%
\special{pa 3339 1288}%
\special{pa 3664 874}%
\special{fp}%
\special{sh 1}%
\special{pa 3664 874}%
\special{pa 3607 914}%
\special{pa 3631 916}%
\special{pa 3639 939}%
\special{pa 3664 874}%
\special{fp}%
%
\special{pn 8}%
\special{pa 3344 1681}%
\special{pa 3489 1505}%
\special{fp}%
%
\special{pn 8}%
\special{pa 3535 1464}%
\special{pa 3865 1045}%
\special{fp}%
%
\special{pn 8}%
\special{pa 1365 1572}%
\special{pa 1787 1897}%
\special{fp}%
\special{sh 1}%
\special{pa 1787 1897}%
\special{pa 1746 1840}%
\special{pa 1745 1864}%
\special{pa 1722 1872}%
\special{pa 1787 1897}%
\special{fp}%
%
\special{pn 8}%
\special{pa 1087 1339}%
\special{pa 669 997}%
\special{fp}%
%
\special{pn 8}%
\special{pa 1632 2094}%
\special{pa 1184 1753}%
\special{fp}%
%
\special{pn 8}%
\special{pa 942 1536}%
\special{pa 503 1173}%
\special{fp}%
\special{sh 1}%
\special{pa 503 1173}%
\special{pa 542 1231}%
\special{pa 544 1207}%
\special{pa 567 1200}%
\special{pa 503 1173}%
\special{fp}%
%
\special{pn 8}%
\special{pa 2870 930}%
\special{pa 3958 1862}%
\special{fp}%
\special{sh 1}%
\special{pa 3958 1862}%
\special{pa 3920 1803}%
\special{pa 3917 1827}%
\special{pa 3894 1834}%
\special{pa 3958 1862}%
\special{fp}%
%
\special{pn 8}%
\special{pa 3772 2073}%
\special{pa 2664 1153}%
\special{fp}%
\special{sh 1}%
\special{pa 2664 1153}%
\special{pa 2703 1211}%
\special{pa 2705 1187}%
\special{pa 2728 1180}%
\special{pa 2664 1153}%
\special{fp}%
\end{picture}%
\vskip3mm


\noindent If $K$ and $J$ are pass-move-equivalent and if $K$ and $K'$ are isotopic, 
then we also say that $K'$ and $J$ are pass-move-equivalent. 
\end{defn}

Note the following: 
If a 1-link $K$ is obtained by a 1-link $J$ by a pass-move, 
 there are two cases 
such that $K_+$ and $K_-$ are not isotopic 
and such that $K_+$ and $K_-$ are isotopic.  
Note this fact in the case of other local moves. 
\\

The pass-move on 1-links is a local-move on 1-links. 
The reason why we use the word `local' in the term `local-move' is as follows:    
When we change a 1-link $K$ into a 1-link $J$ by one pass-move in  $B$ as above,  
we make a change only in $B$ and 
that we do not impose any requirement on 
$K$ or $L$ other than the change only in $B$. 
For example, we do not impose any requirement about the following facts.  
If $K$ is a knot, in which order are the four arcs put in $K$?   
If $K$ is a link and not a knot, 
to which component of  $K$ does each of the four arcs belong?
See also \cite[Note 2.2.(1)]{OgasaZAlex} for this reason. 
\\

The crossing-change on 1-links is also a local-move on 1-links. 
As you know, in the 1-link case local-moves are very useful for research. 
Thus it is very natural to consider local-moves on high dimensional links. 
Indeed, in the high-dimensional-link case local-moves are also very useful for research. 
In the high-dimensional-link case, we must begin by considering what kind of local-moves are fruitful to research high dimensional links. 
The ribbon-move on 2-links and  
the $(p,q)$-pass-move on $(p+q-1)$-links, which we review in this section, 
are two of natural ones. In this paper we discuss them. 
Local moves on high dimensional knots were defined  
in \cite{Ogasa98n, Ogasa04,Ogasa09}  
and 
have been  researched in 
\cite{KauffmanOgasa, KauffmanOgasaB, Ogasa98n,  Ogasa02,  Ogasa04, Ogasa07,   Ogasa09, OgasaT3, OgasaIH, OgasaZAlex, Ogasaribbontwo}. \\

We explain the local-move on the 2-links after 
we review the definition of ribbon spherical 2-links.
Let $m$ be a positive integer. 
A spherical 2-link $L= ( K_1,...,K_m ) $ is called a {\it ribbon} spherical 2-link 
if $L$ satisfies the following properties. 
\smallbreak \noindent
(1) There is a self-transverse immersion 
$f:D^3_1\amalg...\amalg D^3_m$
$\rightarrow S^4$ 
such that $f(\partial D^3_i)=K_i$.  

\smallbreak \noindent
(2) The singular point set $C$  $(\subset S^4$) 
of $f$ consists of double points. 
$C$ is a disjoint union of 2-discs $D^2_i (i=1,...,k)$, 
where $k$ be a nonnegative integer. 

\smallbreak \noindent(3) 
Let $j\in\{1,...,k\}$. 
Let $f^{-1}(D^2_j)=D^2_{jB}\amalg D^2_{jS}$. 
The 2-disc $D^2_{jS}$ is 
embedded in the interior 
of a 3-disc component $D^3_\alpha$ 
for an integer $\alpha\in\{1,...,m\}$.  
The circle $\partial D^2_{jB}$ is 
embedded in 
the boundary 
of a 3-disc component $D^3_\beta$ 
for an integer $\beta\in\{1,...,m\}$.  
The 2-disc $D^2_{jB}$ is 
embedded 
in the 3-disc component $D^3_\beta$.   
(Note that there are two cases, $\alpha=\beta$ and $\alpha\neq\beta$.)

\bigbreak
It is well-known that 
ribbon spherical 2-links are changed into the trivial 2-link by a kind of local-move. 
The local-move is what we call the ribbon-move, 
whose definition we review below.  
The ribbon-move on smooth closed oriented 2-dimensional submanifolds contained in $S^4$ 
was defined in \cite[Definition 1.1]{Ogasa04}.

\begin{defn}\label{tora}   Let $K_+$ and $K_-$ be 
(not necessarily connected or spherical) 
smooth closed oriented 2-dimensional submanifolds contained in $S^4$. 
We say that $K_-$ is obtained from $K_+$ by one {\it ribbon-move } 
if there is a 4-ball $B$ embedded 
in $S^4$ with the following properties.

\smallbreak \noindent
(1) 
$K_+$ and $K_-$ differ only in $B$. 

\smallbreak \noindent
(2) 
$B\cap K_+$ is drawn as in Figure \ref{review}.1.(1).    
$B\cap K_-$ is drawn as in Figure \ref{review}.1.(2).      
In Figure \ref{review}.1.(1) (respectively, \ref{review}.1.(2))  
we draw  $B$ as the product of a 3-ball  and the $t$-axis, 
where this 3-ball is drawn as the product of a 2-disc and the interval. 
$B\cap K_+$ (respectively, $B\cap K_-$) is diffeomorphic to (an annulus)$\amalg$(a 2-disc), 
where $\amalg$ denotes the disjoint union.   
We draw $B\cap K_*$ by bold curves and $B-K_*$ by fine curves.

\smallbreak  
We now describe Figures \ref{review}.1.(1) and \ref{review}.1.(2) in more detail. 
$B\cap K_+$ (respectively, $B\cap K_-$), 
which is diffeomorphic to  $D^2\amalg (S^1\times [0,1])$,  
satisfies the following conditions.

\begin{figure}
     \begin{center}

    \hskip3cm
\unitlength 0.1in
\begin{picture}(59.30,36.30)(0.80,-36.90)
%
\special{pn 8}%
\special{ar 2829 253 417 186  0.0000000 6.2831853}%
%
\special{pn 20}%
\special{pa 2941 2954}%
\special{pa 2941 2046}%
\special{fp}%
%
\special{pn 20}%
\special{pa 2703 2969}%
\special{pa 2703 2054}%
\special{fp}%
%
\special{pn 20}%
\special{pa 2710 2991}%
\special{pa 2723 2962}%
\special{pa 2750 2946}%
\special{pa 2781 2936}%
\special{pa 2813 2932}%
\special{pa 2845 2934}%
\special{pa 2876 2939}%
\special{pa 2906 2952}%
\special{pa 2928 2974}%
\special{pa 2930 3005}%
\special{pa 2910 3029}%
\special{pa 2881 3042}%
\special{pa 2849 3049}%
\special{pa 2817 3051}%
\special{pa 2786 3048}%
\special{pa 2755 3039}%
\special{pa 2727 3024}%
\special{pa 2711 2996}%
\special{pa 2710 2991}%
\special{sp}%
%
\special{pn 20}%
\special{pa 2703 283}%
\special{pa 2703 1191}%
\special{fp}%
%
\special{pn 20}%
\special{pa 2941 268}%
\special{pa 2941 1183}%
\special{fp}%
%
\special{pn 20}%
\special{ar 2822 246 111 60  0.0000000 6.2831853}%
%
\special{pn 8}%
\special{ar 2829 3051 417 186  0.0000000 6.2831853}%
%
\special{pn 8}%
\special{pa 3253 261}%
\special{pa 3253 3051}%
\special{fp}%
%
\special{pn 8}%
\special{pa 2412 268}%
\special{pa 2412 3051}%
\special{fp}%
%
\special{pn 8}%
\special{ar 1267 261 416 186  0.0000000 6.2831853}%
%
\special{pn 8}%
\special{ar 1267 3058 416 186  0.0000000 6.2831853}%
%
\special{pn 8}%
\special{pa 1691 268}%
\special{pa 1691 3058}%
\special{fp}%
%
\special{pn 8}%
\special{pa 850 276}%
\special{pa 850 3058}%
\special{fp}%
%
\special{pn 8}%
\special{ar 4600 246 416 186  0.0000000 6.2831853}%
%
\special{pn 8}%
\special{ar 4600 3043 416 186  0.0000000 6.2831853}%
%
\special{pn 8}%
\special{pa 5024 253}%
\special{pa 5024 3043}%
\special{fp}%
%
\special{pn 8}%
\special{pa 4183 261}%
\special{pa 4183 3043}%
\special{fp}%
%
\special{pn 20}%
\special{ar 2822 1198 111 60  0.0000000 6.2831853}%
%
\special{pn 20}%
\special{ar 2822 2032 111 59  0.0000000 6.2831853}%
%
\special{pn 20}%
\special{ar 2836 1555 417 186  0.0000000 6.2831853}%
%
\special{pn 8}%
\special{pa 2941 1206}%
\special{pa 4481 1206}%
\special{dt 0.045}%
\special{pa 4481 1206}%
\special{pa 4480 1206}%
\special{dt 0.045}%
%
\special{pn 8}%
\special{pa 2956 2039}%
\special{pa 4496 2039}%
\special{dt 0.045}%
\special{pa 4496 2039}%
\special{pa 4495 2039}%
\special{dt 0.045}%
%
\special{pn 20}%
\special{ar 4570 1206 112 59  0.0000000 6.2831853}%
%
\special{pn 20}%
\special{ar 4570 2046 112 60  0.0000000 6.2831853}%
%
\special{pn 20}%
\special{pa 4451 1213}%
\special{pa 4451 2024}%
\special{fp}%
%
\special{pn 20}%
\special{pa 4696 1236}%
\special{pa 4696 2046}%
\special{fp}%
\put(11.1000,-34.3800){\makebox(0,0)[lb]{t=-0.5}}%
\put(25.9800,-34.3800){\makebox(0,0)[lb]{t=0}}%
\put(43.8400,-34.3800){\makebox(0,0)[lb]{t=0.5}}%
%
\special{pn 8}%
\special{ar 2829 253 417 186  0.0000000 6.2831853}%
%
\special{pn 8}%
\special{pa 2412 268}%
\special{pa 2412 3051}%
\special{fp}%
%
\special{pn 8}%
\special{pa 3253 261}%
\special{pa 3253 3051}%
\special{fp}%
%
\special{pn 8}%
\special{ar 2829 3051 417 186  0.0000000 6.2831853}%
%
\special{pn 20}%
\special{ar 2822 246 111 60  0.0000000 6.2831853}%
%
\special{pn 20}%
\special{pa 2941 268}%
\special{pa 2941 1183}%
\special{fp}%
%
\special{pn 20}%
\special{pa 2703 283}%
\special{pa 2703 1191}%
\special{fp}%
%
\special{pn 20}%
\special{pa 2710 2991}%
\special{pa 2723 2962}%
\special{pa 2750 2946}%
\special{pa 2781 2936}%
\special{pa 2813 2932}%
\special{pa 2845 2934}%
\special{pa 2876 2939}%
\special{pa 2906 2952}%
\special{pa 2928 2974}%
\special{pa 2930 3005}%
\special{pa 2910 3029}%
\special{pa 2881 3042}%
\special{pa 2849 3049}%
\special{pa 2817 3051}%
\special{pa 2786 3048}%
\special{pa 2755 3039}%
\special{pa 2727 3024}%
\special{pa 2711 2996}%
\special{pa 2710 2991}%
\special{sp}%
%
\special{pn 20}%
\special{pa 2703 2969}%
\special{pa 2703 2054}%
\special{fp}%
%
\special{pn 20}%
\special{pa 2941 2954}%
\special{pa 2941 2046}%
\special{fp}%
%
\special{pn 8}%
\special{ar 1267 261 416 186  0.0000000 6.2831853}%
%
\special{pn 8}%
\special{ar 1267 3058 416 186  0.0000000 6.2831853}%
%
\special{pn 8}%
\special{pa 1691 268}%
\special{pa 1691 3058}%
\special{fp}%
%
\special{pn 8}%
\special{pa 850 276}%
\special{pa 850 3058}%
\special{fp}%
%
\special{pn 8}%
\special{ar 4600 246 416 186  0.0000000 6.2831853}%
%
\special{pn 8}%
\special{ar 4600 3043 416 186  0.0000000 6.2831853}%
%
\special{pn 8}%
\special{pa 5024 253}%
\special{pa 5024 3043}%
\special{fp}%
%
\special{pn 8}%
\special{pa 4183 261}%
\special{pa 4183 3043}%
\special{fp}%
\put(23.1000,-38.3000){\makebox(0,0)[lb]{$K_+\cap B$}}%
\put(23.1000,-38.3000){\makebox(0,0)[lb]{$K_+\cap B$}}%
%
\special{pn 8}%
\special{pa 80 3540}%
\special{pa 6010 3540}%
\special{fp}%
\special{sh 1}%
\special{pa 6010 3540}%
\special{pa 5943 3520}%
\special{pa 5957 3540}%
\special{pa 5943 3560}%
\special{pa 6010 3540}%
\special{fp}%
\put(54.6000,-38.6000){\makebox(0,0)[lb]{The $t$-axis}}%
\end{picture}%

\vskip4mm \hskip20mm
{\bf Figure \ref{review}.1.(1): Ribbon-move}
\vskip5mm
     \hskip3cm
\unitlength 0.1in
\begin{picture}(60.60,37.30)(1.70,-37.90)
%
\special{pn 8}%
\special{ar 2829 253 417 186  0.0000000 6.2831853}%
%
\special{pn 20}%
\special{pa 2941 2954}%
\special{pa 2941 2046}%
\special{fp}%
%
\special{pn 20}%
\special{pa 2703 2969}%
\special{pa 2703 2054}%
\special{fp}%
%
\special{pn 20}%
\special{pa 2710 2991}%
\special{pa 2723 2962}%
\special{pa 2750 2946}%
\special{pa 2781 2936}%
\special{pa 2813 2932}%
\special{pa 2845 2934}%
\special{pa 2876 2939}%
\special{pa 2906 2952}%
\special{pa 2928 2974}%
\special{pa 2930 3005}%
\special{pa 2910 3029}%
\special{pa 2881 3042}%
\special{pa 2849 3049}%
\special{pa 2817 3051}%
\special{pa 2786 3048}%
\special{pa 2755 3039}%
\special{pa 2727 3024}%
\special{pa 2711 2996}%
\special{pa 2710 2991}%
\special{sp}%
%
\special{pn 20}%
\special{pa 2703 283}%
\special{pa 2703 1191}%
\special{fp}%
%
\special{pn 20}%
\special{pa 2941 268}%
\special{pa 2941 1183}%
\special{fp}%
%
\special{pn 20}%
\special{ar 2822 246 111 60  0.0000000 6.2831853}%
%
\special{pn 8}%
\special{ar 2829 3051 417 186  0.0000000 6.2831853}%
%
\special{pn 8}%
\special{pa 3253 261}%
\special{pa 3253 3051}%
\special{fp}%
%
\special{pn 8}%
\special{pa 2412 268}%
\special{pa 2412 3051}%
\special{fp}%
%
\special{pn 8}%
\special{ar 1267 261 416 186  0.0000000 6.2831853}%
%
\special{pn 8}%
\special{ar 1267 3058 416 186  0.0000000 6.2831853}%
%
\special{pn 8}%
\special{pa 1691 268}%
\special{pa 1691 3058}%
\special{fp}%
%
\special{pn 8}%
\special{pa 850 276}%
\special{pa 850 3058}%
\special{fp}%
%
\special{pn 8}%
\special{ar 4600 246 416 186  0.0000000 6.2831853}%
%
\special{pn 8}%
\special{ar 4600 3043 416 186  0.0000000 6.2831853}%
%
\special{pn 8}%
\special{pa 5024 253}%
\special{pa 5024 3043}%
\special{fp}%
%
\special{pn 8}%
\special{pa 4183 261}%
\special{pa 4183 3043}%
\special{fp}%
%
\special{pn 20}%
\special{ar 2822 1198 111 60  0.0000000 6.2831853}%
%
\special{pn 20}%
\special{ar 2822 2032 111 59  0.0000000 6.2831853}%
%
\special{pn 20}%
\special{ar 2836 1555 417 186  0.0000000 6.2831853}%
\put(11.1000,-34.3800){\makebox(0,0)[lb]{t=-0.5}}%
\put(25.9800,-34.3800){\makebox(0,0)[lb]{t=0}}%
\put(43.8400,-34.3800){\makebox(0,0)[lb]{t=0.5}}%
%
\special{pn 8}%
\special{ar 2829 253 417 186  0.0000000 6.2831853}%
%
\special{pn 8}%
\special{pa 2412 268}%
\special{pa 2412 3051}%
\special{fp}%
%
\special{pn 8}%
\special{pa 3253 261}%
\special{pa 3253 3051}%
\special{fp}%
%
\special{pn 8}%
\special{ar 2829 3051 417 186  0.0000000 6.2831853}%
%
\special{pn 20}%
\special{ar 2822 246 111 60  0.0000000 6.2831853}%
%
\special{pn 20}%
\special{pa 2941 268}%
\special{pa 2941 1183}%
\special{fp}%
%
\special{pn 20}%
\special{pa 2703 283}%
\special{pa 2703 1191}%
\special{fp}%
%
\special{pn 20}%
\special{pa 2710 2991}%
\special{pa 2723 2962}%
\special{pa 2750 2946}%
\special{pa 2781 2936}%
\special{pa 2813 2932}%
\special{pa 2845 2934}%
\special{pa 2876 2939}%
\special{pa 2906 2952}%
\special{pa 2928 2974}%
\special{pa 2930 3005}%
\special{pa 2910 3029}%
\special{pa 2881 3042}%
\special{pa 2849 3049}%
\special{pa 2817 3051}%
\special{pa 2786 3048}%
\special{pa 2755 3039}%
\special{pa 2727 3024}%
\special{pa 2711 2996}%
\special{pa 2710 2991}%
\special{sp}%
%
\special{pn 20}%
\special{pa 2703 2969}%
\special{pa 2703 2054}%
\special{fp}%
%
\special{pn 20}%
\special{pa 2941 2954}%
\special{pa 2941 2046}%
\special{fp}%
%
\special{pn 8}%
\special{ar 1267 261 416 186  0.0000000 6.2831853}%
%
\special{pn 8}%
\special{ar 1267 3058 416 186  0.0000000 6.2831853}%
%
\special{pn 8}%
\special{pa 1691 268}%
\special{pa 1691 3058}%
\special{fp}%
%
\special{pn 8}%
\special{pa 850 276}%
\special{pa 850 3058}%
\special{fp}%
%
\special{pn 8}%
\special{ar 4600 246 416 186  0.0000000 6.2831853}%
%
\special{pn 8}%
\special{ar 4600 3043 416 186  0.0000000 6.2831853}%
%
\special{pn 8}%
\special{pa 5024 253}%
\special{pa 5024 3043}%
\special{fp}%
%
\special{pn 8}%
\special{pa 4183 261}%
\special{pa 4183 3043}%
\special{fp}%
%
\special{pn 8}%
\special{pa 2688 1198}%
\special{pa 1401 1198}%
\special{dt 0.045}%
\special{pa 1401 1198}%
\special{pa 1402 1198}%
\special{dt 0.045}%
%
\special{pn 8}%
\special{pa 2695 2054}%
\special{pa 1423 2054}%
\special{dt 0.045}%
\special{pa 1423 2054}%
\special{pa 1424 2054}%
\special{dt 0.045}%
%
\special{pn 20}%
\special{ar 1274 1198 112 60  0.0000000 6.2831853}%
%
\special{pn 20}%
\special{ar 1274 2039 112 60  0.0000000 6.2831853}%
%
\special{pn 20}%
\special{pa 1155 1206}%
\special{pa 1155 2017}%
\special{fp}%
%
\special{pn 20}%
\special{pa 1401 1228}%
\special{pa 1401 2039}%
\special{fp}%
\put(24.8000,-39.6000){\makebox(0,0)[lb]{$K_-\cap B$}}%
%
\special{pn 8}%
\special{pa 170 3570}%
\special{pa 6230 3560}%
\special{fp}%
\special{sh 1}%
\special{pa 6230 3560}%
\special{pa 6163 3540}%
\special{pa 6177 3560}%
\special{pa 6163 3580}%
\special{pa 6230 3560}%
\special{fp}%
\put(52.8000,-38.9000){\makebox(0,0)[lb]{The $t$-axis}}%
\end{picture}%

\vskip3mm\hskip20mm
{\bf Figure \ref{review}.1.(2): Ribbon-move}
\end{center}
\end{figure}

\begin{itemize}  \item[{\color{white}*}]
We regard $B$ as 
(a closed 2-disc)$\times[0,1]\times\{t| -1\leqq t\leqq1\}$.
Let $B_t\newline
=$(a closed 2-disc)$\times[0,1]\times\{t \}$.  
Note that $B=\cup B_t$. 
In Figure \ref{review}.1.(1) (respectively, \ref{review}.1.(2)),    
we draw 
$B_{-0.5}$ with $B_{-0.5}\cap K_+$, 
$B_{0}$ with $B_{0}\cap K_+$, 
and 
$B_{0.5}$ with $B_{0.5}\cap K_+$  \newline 
(respectively, 
$B_{-0.5}$ with $B_{-0.5}\cap K_-$, 
$B_{0}$ with $B_{0}\cap K_-$, 
and 
$B_{0.5}$ with $B_{0.5}\cap K_-$).   \newline 
We draw $B_{*}\cap K_+$ and $B_{\#}\cap K_-$ by the bold line, 
where  $*,\#\in\{0.5, 0, -0.5\}$. 
We draw $\partial B_t$ by the fine line.

\smallbreak
$B\cap K_+$ has the following properties:  
$B_t\cap K_+$ is empty for $-1\leqq t<0$ and $0.5<t\leqq1$.
$B_0\cap K_+$ is diffeomorphic to 
$D^2\amalg(S^1\times [0,0.3])\amalg(S^1\times [0.7,1])$. 
$B_{0.5}\cap K_+$ is diffeomorphic to $(S^1\times [0.3,0.7])$. 
$B_t\cap K_+$ is diffeomorphic to $S^1\amalg S^1$ for $0<t<0.5$. 
(Here we draw $S^1\times [0,1]$ to have the corner 
in $B_0$ and in $B_{0.5}$. 
However we can let $B\cap K_+$ in $B$ be a smooth submanifold  
by making the corner smooth naturally.)

\smallbreak
$B\cap K_-$ has the following properties:  
$B_t\cap  K_-$ is empty for $-1\leqq t<-0.5$ and $0<t\leqq1$.
$B_0\cap K_-$ is diffeomorphic to 
$D^2\amalg(S^1\times [0, 0.3])\amalg(S^1\times [0.7, 1])$. 
$B_{-0.5}\cap  K_-$ is diffeomorphic to $(S^1\times [0.3, 0.7])$. 
$B_t\cap  K_-$ is diffeomorphic to $S^1\amalg S^1$ for $-0.5<t<0$. 

\smallbreak
In Figure \ref{review}.1.(1) (respectively, \ref{review}.1.(2))  
there are an oriented cylinder $S^1\times [0,1]$ 
and an oriented disc $D^2$ as we stated above. 
We do not make any assumption about 
the orientation of the cylinder and the disc. 
(Of course it holds that 
this orientation of 
(the cylinder)$\amalg$(the disc) 
coincides with 
the orientation of $B\cap K_+$ (respectively, \newline$B\cap K_-$ ).)  
\end{itemize}

Suppose that $K_-$ is obtained from $K_+$ by one ribbon-move 
and that $K'_-$ is isotopic to $K_-$.   
Then we also say that $K'_-$ is obtained from $K_+$ 
by one {\it ribbon-move}.   
If $K_+$ is obtained from $K_-$ by one ribbon-move,  
then we also say that $K_-$ is obtained from $K_+$ by one {\it ribbon-move}.   
$K_+$ and $K_-$ are said to be {\it ribbon-move equivalent} 
if there are 2-knots \newline
$K_+=\bar{K}_1, \bar{K}_2,...,\bar{K}_{r-1},\bar{K}_r=K_-$, 
where $r$ is a positive integer,   
such that 
$\bar{K}_i$ is obtained from $\bar{K}_{i-1}$ $(1< i\leqq r)$ by one ribbon-move. 
\end{defn}


%
%




\bigbreak 
We have the following. 

\begin{thm}\label{kyo}   
$(1)$ 
{\bf(\cite[Main Theorem B]{Ogasa04}.)}  
Not all spherical 2-knots are ribbon-move-equivalent to the trivial 2-knot. 

\smallbreak \noindent
$(2)$  
{\bf(\cite[Note (2) to Main Theorem B]{Ogasa04}.)}  
There is a nonribbon spherical 2-knot which is ribbon-move-equivalent to the trivial 2-knot. 
\end{thm}

Thus it is very natural to consider the following problem.   

\begin{prob}\label{asu}
Classify 2-links (respectively, 2-knots) up to ribbon-move equivalence. 
\end{prob}

 \noindent
In \cite{Ogasa04, Ogasa07, OgasaT3} many partial solutions to Problem \ref{asu} are obtained. 
However this problem remains open.

\bigbreak
The $(p, q)$-pass-move on $n$-knots ($p$ and $q$ are positive integers and  $p+q=n+1$) 
was defined in \cite[section 3]{Ogasa98n},  
and has been studied in 
\cite{KauffmanOgasa, Ogasa98n,  Ogasa02,  Ogasa04, Ogasa07,   Ogasa09, OgasaT3, OgasaIH}. 
We review the definition of the $(p, q)$-pass-move below.    
The (1,1)-pass-move on 1-links is the pass-move on 1-links  
(the definition of the pass-move on 1-links is reviewed in \S\ref{amain}). 
We review the relation between the ribbon-move on 2-links and the (1,2)-pass-move on 2-links 
 (Theorem \ref{r12}).

Each of Figures \ref{review}.2 and \ref{review}.3, 
which consists of the two figures (1) (2), 
represents a diagram of the $(p, q)$-pass-move,  where $q=n+1-p$.  
The figures are generalization of the figure of 
the pass-move on 1-links which is drawn in \S\ref{amain}. 
In Definition \ref{lion} 
we explain the definition of the $(p, q)$-pass-move 
in more detail.

\bigbreak
We use the terms `handle' and `surgeries' in this paper. 
See \cite{Browder, Luck, Smale, Wall} 
for the definition of handles (respectively, surgeries, the attaching parts of handles, 
the attached part, 
other related terms to handles).  
Note that 
an $a$-dimensional $q$-handle $h^q$ is diffeomorphic to $B^q\x B^{a-q}$ (respectively, $B^a$), 
where $B^r$ denotes the $r$-ball, 
and that 
the attaching part of $h^q$ is diffeomorphic $S^{q-1}\x B^{a-q}$. 
Suppose that 
$W$ is obtained by attaching a $k$-handle $h^k$ to $W_0$ 
and that 
$W=W_0\cup h^k$ is embedded in a manifold $M$. 
We say that the submanifold $W\subset M$ 
is obtained by attaching an {\it embedded $k$-handle} $h^k\subset M$ 
to the submanifold $W_0\subset M$.

\begin{defn}\label{lion}  
Let $n$ and $p$ be positive integers. 
Let $n+1-p>0.$     
We now define a $(p, n+1-p)$-pass-move in a $(n+2)$-ball. We first explain 
Figures \ref{review}.2 and \ref{review}.3. 

Regard an $(n+2)$-ball $B=D^{n+2}$ as $[-1,1]\x D^p\x D^{n+1-p}$ 
as drawn in Figure \ref{review}.3.  

Attach an embedded $(n+1)$-dimensional $(n+1-p)$-handle 
$h^{n+1-p} \subset\{0\}\x D^p\x D^{n+1-p}$ 
to 
$\{0\}\x\partial(D^p\x D^{n+1-p})$
along 
$\{0\}\x\{*\}\x \partial D^{n+1-p}$
so that 
the core of $h^{n+1-p}$ coincides with 
$\{0\}\x\{*\}\x D^{n+1-p}$ (see Figures \ref{review}.2 and \ref{review}.3).

Attach an embedded $(n+1)$-dimensional $p$-handle 
$h^p \subset\{0\}\x D^p\x D^{n+1-p}$ 
to \newline
$\{0\}\x\partial(D^p\x D^{n+1-p})$
along 
$\{0\}\x \partial D^p\x\{*\}$
so that 
the core of $h^p$ coincides with 
$\{0\}\x D^p\x\{*\}$.  
Note $h^{n+1-p}\cap h^p\neq\phi.$

Move $h^p$ in $B$ by using an isotopy with keeping $h^p\cap \partial B$, 
 let 
the resultant submanifold contained in $\{t\geqq0\}\x D^p\x D^{n+1-p}$ 
(respectively, $\{t\leqq0\}\x D^p\x D^{n+1-p}$), 
and call the submanifold $h_+^p$ (respectively, $h_-^p$).  
Suppose that $h_+^p\cap h^{n+1-p}=\phi$ and 
that $h_-^p\cap h^{n+1-p}=\phi$.   
(see Figures \ref{review}.2 and \ref{review}.3).

%
%

Let $K_+$ and $K_-$ be 
$n$-dimensional closed oriented submanifolds contained in $S^{n+2}$. 
Embed the $(n+2)$-ball $B$ 
in $S^{n+2}$. 
Let $K_+$ and $K_-$ differ only in $B$. 
Let $K_+$ (respectively, $K_-$) satisfy the condition  \newline
\hskip39mm$K_+\cap \mathrm{Int}B=
(\partial h^{p}_+-\partial B)\cup(\partial h^{n+1-p}-\partial B)$ \newline  
\hskip28mm 
${\rm(respectively,}\hskip1mm K_-\cap \mathrm{Int}B=
(\partial h^{p}_--\partial B)\cup(\partial h^{n+1-p}-\partial B)),$  \newline
where we suppose that there is not  
$h^{p}_-$  
(respectively, $h^{p}_+$)   
in $B$.  
Then we say that 
$K_+$ (respectively, $K_-$)  is obtained from $K_-$ (respectively, $K_+$) 
by one {\it $(p,n+1-p)$-pass-move} in $B$.   

\end{defn}

In Definition \ref{lion}, we have the following: 
Let $\sharp\in\{+,-\}.$   
There is a Seifert hypersurface $V_\sharp\subset S^{n+2}$ for $K_\sharp$ such that 
$V_\sharp\cap B=h_\sharp^p\cup h^{n+1-p}.$ 
\noindent 
(The idea of the proof is Thom-Pontrjagin construction.)
We say that  
$V_-$ (respectively, $V_+$) is obtained from $V_+$ (respectively, $V_-$) 
by a {\it $(p,n+1-p)$-pass-move} in $B$.

\bigbreak
In Definition \ref{lion}, note the following: 
Let 
$V_0=V_\sharp-\text{Int}B$\newline  
$=\text{the closure of }`V_\sharp- (h_\sharp^p\cup h^{n+1-p})'\text{ in }{S^{n+2}}.$  \\
We can say that 
we attach 
an embedded $(n+1)$-dimensional $p$-handle $h_\#^p\subset S^{n+2}$ 
and 
an embedded $(n+1)$-dimensional $(n+1-p)$-handle $h^{n+1-p}\subset S^{n+2}$ 
to the $(n+1)$-dimensional submanifold $V_0\subset S^{n+2}$, 
and 
obtain the $(n+1)$-dimensional submanifold $V_\#\subset S^{n+2}$.

\begin{figure}
\unitlength 0.1in
\begin{picture}(52.00,24.84)(1.50,-32.30)
%
\special{pn 8}%
\special{ar 4488 1824 843 843  0.8934700 6.2831853}%
\special{ar 4488 1824 843 843  0.0000000 0.8617842}%
%
\special{pn 8}%
\special{pa 3829 1312}%
\special{pa 4891 2556}%
\special{fp}%
%
\special{pn 8}%
\special{pa 5155 2326}%
\special{pa 4104 1083}%
\special{fp}%
%
\special{pn 8}%
\special{pa 4891 1083}%
\special{pa 4515 1440}%
\special{fp}%
\special{pa 4168 1806}%
\special{pa 3738 2208}%
\special{fp}%
\special{pa 5147 1320}%
\special{pa 4735 1686}%
\special{fp}%
\special{pa 4369 2071}%
\special{pa 3966 2455}%
\special{fp}%
%
\special{pn 8}%
\special{ar 1652 1858 843 843  5.5899622 6.2831853}%
\special{ar 1652 1858 843 843  0.0000000 5.5598920}%
%
\special{pn 8}%
\special{pa 904 1468}%
\special{pa 1268 1838}%
\special{fp}%
\special{pa 1639 2180}%
\special{pa 2047 2604}%
\special{fp}%
\special{pa 1138 1209}%
\special{pa 1512 1614}%
\special{fp}%
\special{pa 1901 1974}%
\special{pa 2291 2371}%
\special{fp}%
%
\special{pn 8}%
\special{pa 2146 1183}%
\special{pa 917 2254}%
\special{fp}%
%
\special{pn 8}%
\special{pa 1149 2526}%
\special{pa 2377 1445}%
\special{fp}%
%
\special{pn 8}%
\special{pa 884 3002}%
\special{pa 1039 2179}%
\special{dt 0.045}%
\special{sh 1}%
\special{pa 1039 2179}%
\special{pa 1007 2241}%
\special{pa 1029 2231}%
\special{pa 1046 2248}%
\special{pa 1039 2179}%
\special{fp}%
\special{pa 921 2983}%
\special{pa 1258 2426}%
\special{dt 0.045}%
\special{sh 1}%
\special{pa 1258 2426}%
\special{pa 1206 2473}%
\special{pa 1230 2472}%
\special{pa 1241 2493}%
\special{pa 1258 2426}%
\special{fp}%
%
\special{pn 8}%
\special{pa 2402 2801}%
\special{pa 2228 2334}%
\special{dt 0.045}%
\special{sh 1}%
\special{pa 2228 2334}%
\special{pa 2233 2403}%
\special{pa 2247 2384}%
\special{pa 2270 2389}%
\special{pa 2228 2334}%
\special{fp}%
\special{pa 2393 2783}%
\special{pa 1853 2361}%
\special{dt 0.045}%
\special{sh 1}%
\special{pa 1853 2361}%
\special{pa 1893 2418}%
\special{pa 1895 2394}%
\special{pa 1918 2386}%
\special{pa 1853 2361}%
\special{fp}%
\put(22.1900,-30.6600){\makebox(0,0)[lb]{$S^{p-1}\x D^{n+1-p}$}}%
\put(3.0700,-31.6600){\makebox(0,0)[lb]{$D^p\x S^{n-p}$}}%
\put(21.8200,-33.0500){\makebox(0,0)[lb]{$=\overline{\partial h^{n+1-p}-\partial B}$}}%
\put(11.9400,-9.5300){\makebox(0,0)[lb]{$B\cap K_+$}}%
\put(40.6700,-9.1600){\makebox(0,0)[lb]{$B\cap K_-$}}%
\put(1.5000,-34.0000){\makebox(0,0)[lb]{$=\overline{\partial h^p_+-\partial B}$}}%
%
\special{pn 8}%
\special{pa 1100 1450}%
\special{pa 1111 1419}%
\special{pa 1121 1389}%
\special{pa 1128 1358}%
\special{pa 1130 1326}%
\special{pa 1128 1295}%
\special{pa 1122 1263}%
\special{pa 1112 1232}%
\special{pa 1100 1201}%
\special{pa 1087 1171}%
\special{pa 1071 1143}%
\special{pa 1053 1116}%
\special{pa 1033 1092}%
\special{pa 1009 1072}%
\special{pa 982 1055}%
\special{pa 954 1040}%
\special{pa 924 1026}%
\special{pa 893 1015}%
\special{pa 862 1004}%
\special{pa 831 1000}%
\special{pa 800 1007}%
\special{pa 770 1020}%
\special{sp}%
\put(4.1000,-11.9000){\makebox(0,0)[lb]{$h^{n+1-p}$}}%
\put(32.6000,-11.8000){\makebox(0,0)[lb]{$h^{n+1-p}$}}%
%
\special{pn 8}%
\special{pa 4040 1360}%
\special{pa 4035 1327}%
\special{pa 4028 1296}%
\special{pa 4018 1266}%
\special{pa 4002 1238}%
\special{pa 3982 1214}%
\special{pa 3961 1189}%
\special{pa 3941 1165}%
\special{pa 3921 1140}%
\special{pa 3902 1115}%
\special{pa 3883 1089}%
\special{pa 3864 1062}%
\special{pa 3844 1035}%
\special{pa 3821 1012}%
\special{pa 3794 998}%
\special{pa 3761 994}%
\special{pa 3729 1000}%
\special{pa 3700 1015}%
\special{pa 3670 1020}%
\special{sp}%
%
\special{pn 8}%
\special{pa 2190 1380}%
\special{pa 2199 1349}%
\special{pa 2209 1318}%
\special{pa 2219 1288}%
\special{pa 2231 1259}%
\special{pa 2245 1230}%
\special{pa 2261 1202}%
\special{pa 2280 1176}%
\special{pa 2300 1151}%
\special{pa 2322 1127}%
\special{pa 2346 1106}%
\special{pa 2373 1089}%
\special{pa 2404 1078}%
\special{pa 2435 1069}%
\special{pa 2464 1058}%
\special{pa 2488 1040}%
\special{pa 2508 1016}%
\special{pa 2524 987}%
\special{pa 2539 956}%
\special{pa 2550 930}%
\special{sp}%
\put(25.0000,-9.5000){\makebox(0,0)[lb]{$h^p_+$}}%
\put(53.5000,-12.4000){\makebox(0,0)[lb]{$h^p_-$}}%
%
\special{pn 8}%
\special{pa 5000 1290}%
\special{pa 5028 1274}%
\special{pa 5054 1255}%
\special{pa 5078 1234}%
\special{pa 5098 1208}%
\special{pa 5117 1180}%
\special{pa 5140 1160}%
\special{pa 5170 1152}%
\special{pa 5203 1145}%
\special{pa 5234 1140}%
\special{pa 5265 1145}%
\special{pa 5297 1150}%
\special{pa 5320 1150}%
\special{sp}%
\end{picture}%

\vskip5mm
\hskip5mm\text{{\bf  Figure \ref{review}.2: 
The $(p, n+1-p)$-pass-move on an $n$-dimensional closed submanifold  \newline
 }}

\hskip-10mm\text{ {\bf 
contained in $S^{n+2}$. 
Note $B=B^{n+2}=D^{n+2}\subset S^{n+2}$.   
}}

\end{figure}

\begin{figure}
\begin{center}
\includegraphics[width=107mm]{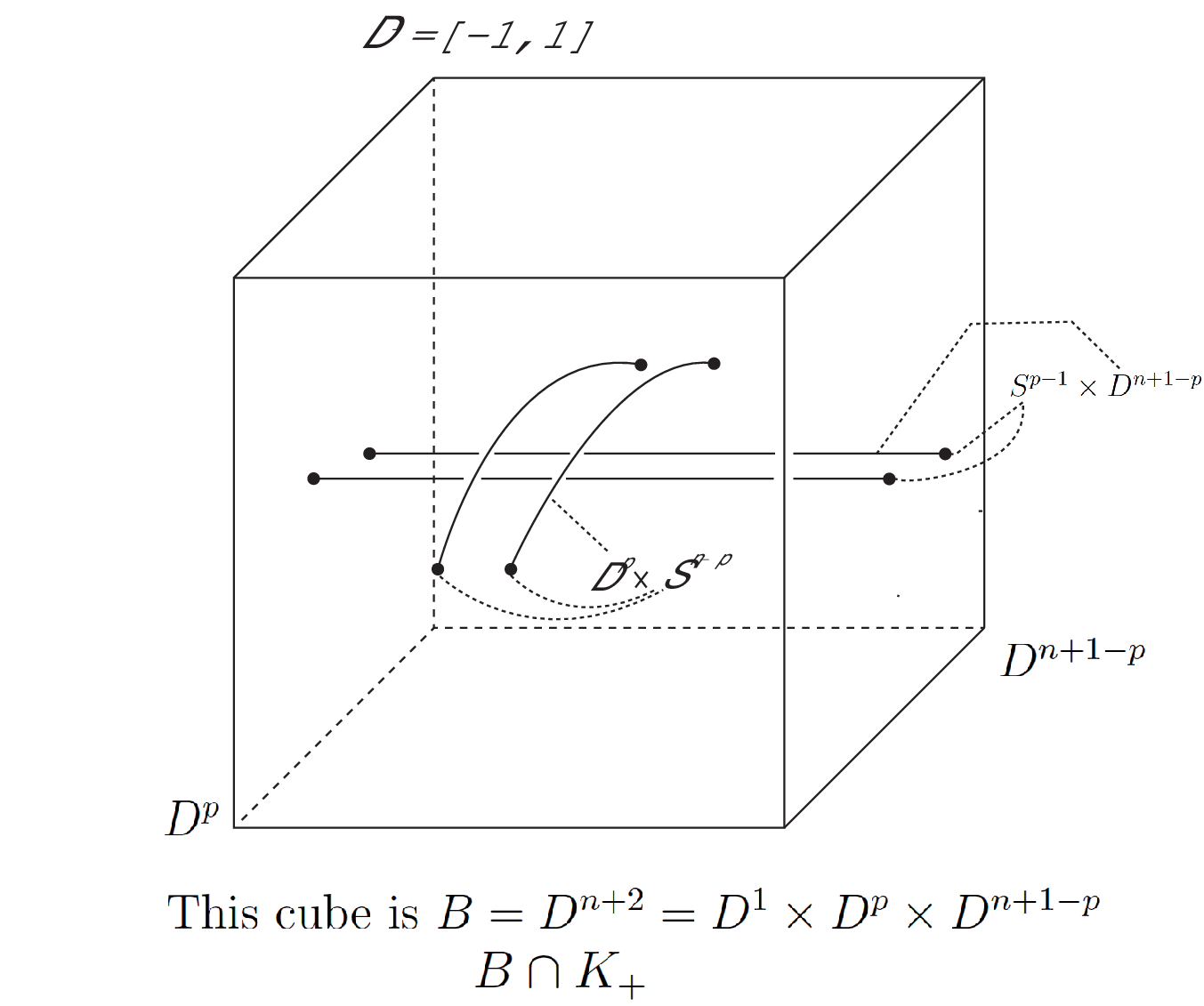}

\text{
Figure \ref{review}.3.(1): {\bf The $(p,n+1-p)$-pass-move}
}

\vskip10mm
\includegraphics[width=70mm]{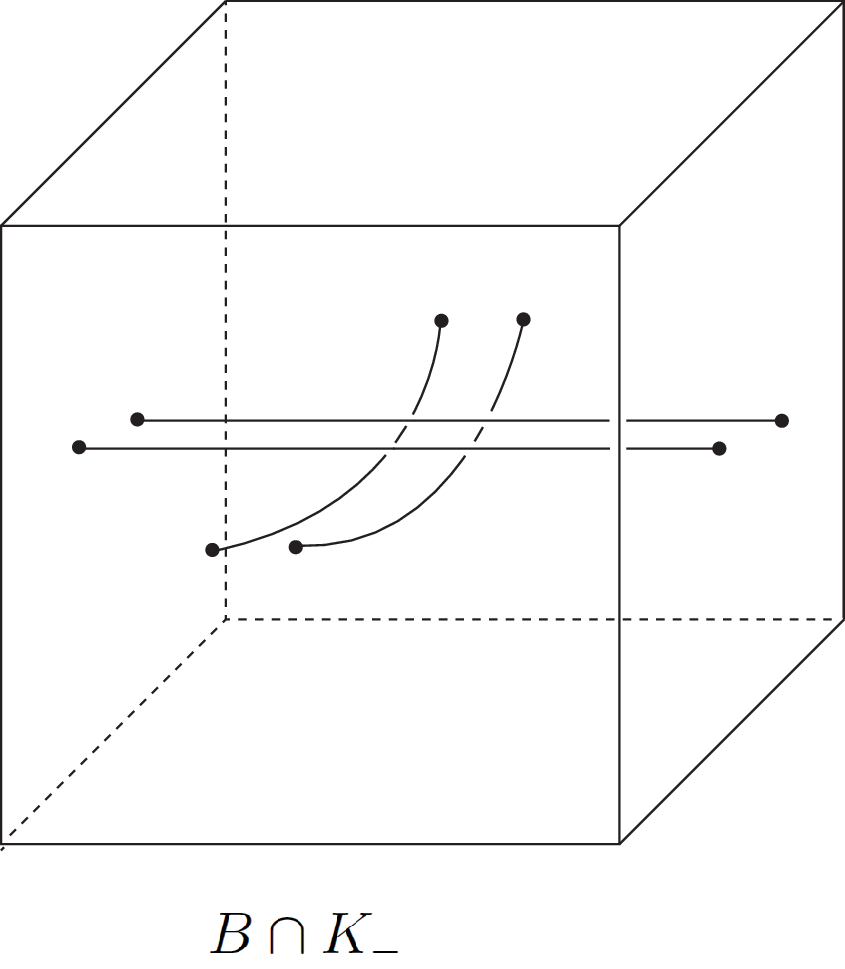}

\vskip2mm
\text{Figure \ref{review}.3.(2): {\bf The $(p,n+1-p)$-pass-move} }
\end{center}
\end{figure}

We have the following. 

\begin{thm}\label{r12} {\rm{\bf(\cite[Propositions 4.2 and 4.3]{Ogasa04}.)}}
  Let $K$ and $K'$ be 2-dimensional closed oriented submanifolds contained in $S^4$.  

\smallbreak \noindent 
$(1)$  The following conditions {\rm (i)} and {\rm (ii)} are equivalent. 

\smallbreak \noindent{\color{white}$(1)$}   
${\rm (i)}$  
$K$ is $($1,2$)$-pass-move-equivalent to $K'$.

\smallbreak \noindent{\color{white}$(1)$}   
${\rm (ii)}$    
$K$ is  ribbon-move-equivalent to $K'$. 

\smallbreak \noindent $(2)$ 
If $K$ is obtained from $K'$ by one ribbon-move, 
then $K$ is obtained from $K'$ by one 

\smallbreak \noindent{\color{white}$(2)$}   
$($1,2$)$-pass-move.
\end{thm}

We draw the (1,2)-pass-move (the case $p=1$, $q=2$, and $n=2$) 
in Figure \ref{review}.4, which consists of the two figures (1) (2), as follows:    
Let $K_+$ be obtained from $K_-$ by a $(1,2)$-pass-move in $B$. 
Thus $B$ has the following properties according to Definition \ref{lion}. 

\smallbreak
\noindent 
(i) $K_+$ and $K_-$differ only in $B$.   

\smallbreak
\noindent 
(ii) $B\cap K_+$ and $B\cap K_-$ are shown in Figure \ref{review}.4.   

\begin{itemize}\item[]
We regard $B$ as 
(a 2-disc)$\x[0,1]\x\{t\vert -1\leqq t\leqq1\}$.  
Let $*=+,-$. 
We draw 
$B_{-0.5}\cap K_*$, $B_0\cap K_*$, $B_{0.5}\cap K_*$, 
where $B_\xi$=(a 2-disc)$\x[0,1]\x\{t\vert t=\xi\}$.  
We suppose that 
each vector $\overrightarrow{x}$, $\overrightarrow{y}$ 
in Figure \ref{review}.4  
is a tangent vector of each disc at a point. 
(Note that we use $\overrightarrow{x}$ (respectively, $\overrightarrow{y}$)
for different vectors.)
The orientation of each disc 
in Figure \ref{review}.4  
is determined by the each ordered set 
$\{\overrightarrow{x},\overrightarrow{y}\}$. 
See more explanation near \cite[Figure 4.1 and 4.2]{Ogasa04}.  
%
\end{itemize}

\begin{figure}
     \begin{center}
 \input 3+.tex  

\vskip5mm\hskip16mm  
{\bf Figure \ref{review}.4.(1): The (1,2)-pass-move} 
\vskip7mm

 \input 3-.tex  

\vskip5mm\hskip16mm  
{\bf Figure \ref{review}.4.(2): The (1,2)-pass-move}

     \end{center}
\end{figure}




\bigbreak
Local moves on high dimensional submanifolds are exciting ways of explicit construction of high dimensional figures. They are also a generalization of local-moves on 1-links as we saw in this section. 
They are useful to research link cobordism, knot cobordism, and the intersection of submanifolds (see \cite{Ogasa98n}).   
There remain many exciting problems. Some of them are proper in high dimension and others are analogous to 1-dimensional cases. For example, we do not know a local-move on high dimensional knots which is an unknotting operation.

\bigbreak
We proved the following relations between local moves on high dimensional knots and several invariants of high dimensional knots. 

In \cite[Theorem 4.1]{Ogasa98n}  the following is proved: 
The $(p+1,p+1)$-pass-move on spherical $(2p+1)$-knots 
preserves the Arf invariant (respectively, the signature) if $p$ is even (respectively, odd). 
Furthermore the following is proved: Let $p$ be even (respectively, odd). 
Simple $(2p+1)$-knots, $K$ and $J$, are $(p+1,p+1)$-pass-equivalent 
if and only if 
their Arf invariant (respectively, their signature) are the same.

\cite[Theorem 4.1]{Ogasa98n} 
and 
\cite[Proposition 12.1]{KauffmanOgasa} 
imply the following:  Let $p$ be even (respectively, odd).
A $(2p+1)$ knot $K$ is $(p+1,p+1)$-pass-move equivalent to the unknot 
if and only if 
$K$ is a simple spherical knot and 
the Arf invariant (resp. the signature) is zero in the case where 
$p$ is even (resp. odd). 
\\

In \cite[Their main theorems]{Ogasa04, Ogasa07,   OgasaT3} the following are proved: 
The (1,2)-pass move (respectively, the ribbon-move) on spherical 2-knots preserves 
the $\mu$-invariant of 2-knots, 
the $\Q/\Z$-valued $\widetilde\eta$-invariants of 2-knots, 
the Farber-Levine pairing of 2-knots,   
and 
partial information of the cup product of three elements in 
$H^1_{\rm cpt}($the complement of each 2-knot$)$.   
Theorem \ref{kyo} 
and the results written below Problem \ref{asu} 
were proved by using a few of these results.
\\

In \cite[Theorems 3.2 and 3.3]{Ogasa09} the following is proved:  
For the Alexander polynomial $A(\hskip2mm)$ of high dimensional knots 
we have an identity 
$$A(K_+)-A(K_-)=(t-1)\cdot A(K_0)$$ 
associated with the twist move on $(4k+1)$-dimensional knots, 
where $k$ is a nonnegative integer 
(respectively, the $(p,q)$-pass-move, where $p\neq q$, on high dimensional knots).   
The twist move is a kind of local moves on high dimensional knots, and 
is defined in \cite[section 3]{Ogasa09}.
\\

In \cite[Theorem 9.2]{KauffmanOgasa} the following is proved: 
For the Alexander polynomial $A(\hskip2mm)$ of 
$(4k+3)$-dimensional knots ($k$ is a nonnegative integer), 
we have an identity 
$$A(K_+)-A(K_-)=(t+1)\cdot A(K_0)$$ 
associated with the twist move.
It is a new type of local move identities of knot polynomial. 
Note $(t+1)$ in the right hand side. It is not $(t-1)$.

\bigbreak
\section{Main results on relations between local  moves on knots and knot products   
}\label{main2}

\noindent
It is very natural to consider the following problem:  
If codimension two closed oriented submanifolds $K$ and $J$ in $S^{l}$ 
are obtained from each other by a local-move, 
then for a codimension two closed oriented submanifold $A$ in $S^{w}$, 
are $K\otimes A$ and  $J\otimes A$ are obtained from each other by another local-move?

In \cite{KauffmanOgasa}  we considered this problem  
and obtained several results associated with this problem.  
It means that we  
combined research of knot products and that of local moves on high dimensional knots.   
Some of the results are as follows. 
In this paper we obtain Main Theorems  \ref{sugar} and \ref{melon}.  
Some other ones are in \cite{KauffmanOgasaB}. \\

Let $k$ be a positive integer. 
Let $P\otimes^kQ$ denote the knot product $P\otimes\underbrace{Q\otimes...\otimes Q}_k$. 
Let  $P\otimes^0Q$ denote $P$. 
In \cite[Theorem 8.1]{KauffmanOgasa} 
we proved the following: 
Let $\mu$ be a positive integer.
If a 1-knot $A$ is obtained from a 1-knot $B$ by one pass-move, 
then 
$A\otimes^{\mu} {\rm Hopf}$ is obtained from $B\otimes^{\mu} {\rm Hopf}$
by one $(2\mu+1, 2\mu+1)$-pass-move.   
\\

In \cite[Theorem 8.10]{KauffmanOgasa} 
we proved the following: 
Let $A$ be a 1-knot. 
Let $\mu$ be a positive integer. 
Let $J=A\otimes^{\mu} {\rm Hopf}$. 
Let $K$ be obtained from $J$ by one $(2\mu+1, 2\mu+1)$-pass-move. 
Then there is a 1-knot $B$ such that  $K=B\otimes^{\mu} {\rm Hopf}$ 
and such that $A$ is pass-equivalent to $B$.  \\

\cite[Theorems 8.1 and 8.10]{KauffmanOgasa} 
imply the following: 
Let $\mu$ be a positive integer.
A 1-knot $A$ is pass-equivalent to a 1-knot $B$ if and only if 
$A\otimes^{\mu} {\rm Hopf}$ is 
$(2\mu+1, 2\mu+1)$-pass-equivalent to 
$B\otimes^{\mu} {\rm Hopf}$.
\\

\cite[Theorem 4.1]{KauffmanOgasa} proved  the following: 
Suppose that two 1-links $J$ and $K$ differ by a single crossing change.
Let $\mu$ be a positive integer. 
Then the knot products, 
$J \otimes^\mu{\rm Hopf}$ 
and   
$K\otimes^\mu{\rm Hopf}$,   
differ by a single twist-move. 
The twist-move is a kind of local moves, which is defined in \cite[section 3]{Ogasa09}, 
as written in the last part of \S\ref{review}. 
\\

 \cite[Theorem 7.1]{KauffmanOgasa} proved the following: 
Let $m$ and $\nu$ be nonnegative integers.  
Suppose that two 
$($not necessarily connected or spherical$)$ $(2m+1)$-dimensional closed oriented 
submanifolds, $J$ and $K$, of $S^{2m+3}$ differ by a single twist-move.
Then the $(2m+2\nu+1)$-submanifolds, 
$J\otimes^\nu[2]$  
and 
$K\otimes^\nu[2]$, 
 of $S^{2m+2\nu+3}$  
differ by a single twist-move.
Note that the twist move on 1-links is the crossing change on those 
(see \cite[section 3]{Ogasa09}).

\begin{note}\label{kouiunomo}
Let $[2]$ be the empty knot of degree two, which is defined in 
\cite[line 7 in section 2]{Kauffman} and 
\cite[the last line of Definition 2.1 and the fifth line of Definition 3.1]{KauffmanNeumann}.   

We have Hopf$=[2]\otimes[2]$.  
Recall that Hopf link is a Brieskorn submanifold as written in \cite{Milnor}. 
See \cite[the sixth line from the bottom of page 1106]{Kauffman} and 
\cite[the last sentence of \S6, and the first remark in page 390]{KauffmanNeumann}. 
\end{note}

 \cite[Theorem 7.3]{KauffmanOgasa} proved the following: 
Let $k$ be a positive integer.    
Let $K$ $($respectively, $J)$ be $(4k+5)$-submanifold of $S^{4k+7}$. 
Suppose that $K$ and $J$ differ by a single twist-move and are not isotopic. 
Suppose that $K$ is isotopic to $A\otimes^{k+1}{\rm Hopf}$ for a 1-knot $A$. 
Then  
there is a unique equivalence class of simple $(4k+1)$-knots for $K$ $($respectively, $J)$   
with the following properties. 

\smallbreak
\noindent $\mathrm{(i)}$
There is a representative element $K'$ of the above equivalence class for $K$ such that 
$K$ is isotopic to $K'\otimes{\rm Hopf}$.   

\smallbreak
\noindent $\mathrm{(ii)}$
There is a representative element $J'$ of the above equivalence class for $J$ such that 
$J$ is isotopic to $J'\otimes{\rm Hopf}$.   

\smallbreak
\noindent $\mathrm{(iii)}$
$K'$ and $J'$ differ by a single twist-move and are not isotopic. 
\\

In this paper we obtain the following results on relations 
between knot products and local moves on knots.

\begin{mthm} \label{sugar}      
Let $J$ and $K$ be 1-links in $S^3$. 
Suppose that $J$ is obtained from $K$ by a single pass-move.  
Let $k$ be a positive integer.   
Then $J\otimes^k{\rm Hopf}$ is obtained from 
$K\otimes^k{\rm Hopf}$  by a single $(2k+1,2k+1)$-pass-move.  
\end{mthm}

\noindent{\bf Note.} 
(1) \cite[Theorem 8.1]{KauffmanOgasa} proved the case where $J$ is a 1-component 1-link. 

\smallbreak \noindent
(2)  
The converse of Main Theorem \ref{sugar} is false by 
\cite[Theorem 1.6.(1)]{Ogasaribbontwo}.
\bigbreak

\noindent
Main Theorem \ref{sugar} implies Theorems \ref{B}.(1) and \ref{B}.(3).   
(We prove Theorem \ref{B} in \S\ref{hightsuzuki}.)

\begin{thm} \label{B}    
$(1)$  Let $a,b,a',b'$ and $k$ be positive integers. 
If the $(a,b)$ torus link is pass-move equivalent to the $(a',b')$ torus link, 
then the Brieskorn manifold $\Sigma(a,b,\underbrace{2,...,2}_{2k})$ 
is diffeomorphic to $\Sigma(a',b',\underbrace{2,...,2}_{2k})$ as submanifolds. 

\smallbreak\noindent$(2)$  
The converse of $(1)$ is false in general.  

\smallbreak\noindent$(3)$  
If two 1-links $J$ and $K$ are pass-move equivalent, 
then $J\otimes^k{\rm Hopf}$ is diffeomorphic to $K\otimes^k{\rm Hopf}$.   
\end{thm}



\begin{mthm} \label{melon}  
Let $J$ and $K$ be $($not necessarily connected or spherical$)$ 2-dimensional closed oriented submanifolds  in $S^4$. Suppose that $J$ is obtained from $K$ by a single ribbon-move.  
Let $k$ be an integer$\geq2$.  
Then 
$J\otimes^k{\rm Hopf}$ 
is obtained from $K\otimes^k{\rm Hopf}$  by a single $(2k+1,2k+2)$-pass-move.  
\end{mthm}

\noindent{\bf Note.} 
The converse of Main Theorem \ref{melon} is false by \cite[Theorem 8.12]{KauffmanOgasa}.
\\

Let $n$ be a positive integer. 
Let $\alpha$ be a sufficiently large positive integer.  
Main Theorems \ref{coffee},  \ref{sugar}, and \ref{melon} connect  
the classification problem of $n$-knots by a local-move-equivalence 
with that of $(n+\alpha)$-dimensional stable knots by another local-move-equivalence.  
Of course the latter classification problem is easier than the former one 
by results in \cite{Farber1978, Farber1980, Farber1983, Farber1984I, Farber1984II}. 
  See also \cite{Farber1978, Farber1980, Farber1983, Farber1984I, Farber1984II} for 
the definition of stable knots.   
(If $n=1,2$, `stable' is replaced by `simple'.)  
Thus we can consider lower and higher dimensional knots together.

\bigbreak
\section{
Proof of Main Theorem \ref{melon}}\label{high}
\noindent   
The idea of the proof is as follows. 
Take a Seifert hypersurface $V_J$ (respectively $V_K$) for $J$ (respectively $K$) 
such that $V_J$ and $V_K$ differ only in the 4-ball where the one ribbon-move is carried out. 
We make a Seifert hypersurface 
$V_{J\otimes^k{\rm Hopf}}$ (respectively $V_{K\otimes^k{\rm Hopf}}$) 
for $J\otimes^k{\rm Hopf}$ (respectively $K\otimes^k{\rm Hopf}$)  
from $V_J$ (respectively $V_K$) 
by using Lemmas \ref{soda} and \ref{sodasui}
such that $V_{J\otimes^k{\rm Hopf}}$ and $V_{K\otimes^k{\rm Hopf}}$ 
differ only by one  $(2k+1,2k+2)$-pass-move in a $(4k+4)$-ball of $S^{4k+4}$.  
Their difference implies that 
$J\otimes^k{\rm Hopf}$ and $K\otimes^k{\rm Hopf}$  
differ by one $(2k+1,2k+2)$-pass-move. 

Let $m$ be a nonnegative integer. 
Recall that Hopf=$[2]\otimes[2]$. 
We point out that there are not 
a Seifert hypersurface $V_{J\otimes^{2m+1}[2]}$ for $J\otimes^{2m+1}[2]$ 
and 
a Seifert hypersurface $V_{K\otimes^{2m+1}[2]}$ for $K\otimes^{2m+1}[2]$ 
which are diffeomorphic 
 in general  
although  the above 
$V_{J\otimes^k{\rm Hopf}}$ 
and 
$V_{K\otimes^k{\rm Hopf}}$ 
have the above property.

\bigbreak
We need the following lemmas in order to prove Main Theorem \ref{melon}.   

\bigbreak
\begin{figure}

\unitlength 0.1in
\begin{picture}(57.55,28.37)(1.50,-36.00)
%
\special{pn 8}%
\special{ar 4968 1942 937 937  0.8932296 6.2831853}%
\special{ar 4968 1942 937 937  0.0000000 0.8617069}%
%
\special{pn 8}%
\special{pa 4236 1373}%
\special{pa 5415 2755}%
\special{fp}%
%
\special{pn 8}%
\special{pa 5709 2500}%
\special{pa 4541 1118}%
\special{fp}%
%
\special{pn 8}%
\special{pa 5415 1118}%
\special{pa 4998 1515}%
\special{fp}%
\special{pa 4612 1922}%
\special{pa 4134 2369}%
\special{fp}%
\special{pa 5700 1382}%
\special{pa 5242 1789}%
\special{fp}%
\special{pa 4835 2216}%
\special{pa 4388 2643}%
\special{fp}%
%
\special{pn 8}%
\special{ar 1817 1980 937 937  5.5898109 6.2831853}%
\special{ar 1817 1980 937 937  0.0000000 5.5597350}%
%
\special{pn 8}%
\special{pa 1258 2722}%
\special{pa 2622 1521}%
\special{fp}%
%
\special{pn 8}%
\special{pa 2365 1230}%
\special{pa 1000 2420}%
\special{fp}%
%
\special{pn 8}%
\special{pa 986 1546}%
\special{pa 1390 1957}%
\special{fp}%
\special{pa 1802 2337}%
\special{pa 2255 2808}%
\special{fp}%
\special{pa 1246 1258}%
\special{pa 1661 1709}%
\special{fp}%
\special{pa 2093 2109}%
\special{pa 2527 2550}%
\special{fp}%
\put(13.0800,-9.7400){\makebox(0,0)[lb]{$B^4\cap J$}}%
\put(45.0000,-9.3300){\makebox(0,0)[lb]{$B^4\cap K$}}%
%
\special{pn 8}%
\special{pa 1160 2250}%
\special{pa 1128 2250}%
\special{pa 1096 2250}%
\special{pa 1064 2249}%
\special{pa 1032 2248}%
\special{pa 1000 2247}%
\special{pa 968 2246}%
\special{pa 936 2244}%
\special{pa 904 2242}%
\special{pa 872 2240}%
\special{pa 840 2238}%
\special{pa 808 2243}%
\special{pa 778 2256}%
\special{pa 752 2276}%
\special{pa 729 2302}%
\special{pa 713 2332}%
\special{pa 704 2363}%
\special{pa 701 2394}%
\special{pa 700 2427}%
\special{pa 698 2459}%
\special{pa 696 2491}%
\special{pa 694 2523}%
\special{pa 691 2554}%
\special{pa 688 2586}%
\special{pa 685 2618}%
\special{pa 682 2650}%
\special{pa 680 2670}%
\special{sp -0.045}%
%
\special{pn 8}%
\special{pa 1400 2620}%
\special{pa 1400 2652}%
\special{pa 1401 2684}%
\special{pa 1402 2716}%
\special{pa 1400 2748}%
\special{pa 1394 2780}%
\special{pa 1382 2810}%
\special{pa 1365 2837}%
\special{pa 1341 2859}%
\special{pa 1316 2879}%
\special{pa 1289 2897}%
\special{pa 1260 2912}%
\special{pa 1231 2924}%
\special{pa 1200 2932}%
\special{pa 1168 2937}%
\special{pa 1136 2939}%
\special{pa 1104 2941}%
\special{pa 1072 2942}%
\special{pa 1040 2942}%
\special{pa 1008 2940}%
\special{pa 976 2936}%
\special{pa 945 2930}%
\special{pa 914 2923}%
\special{pa 900 2920}%
\special{sp -0.045}%
\put(1.5000,-28.9000){\makebox(0,0)[lb]{$S^1\x B^1$}}%
%
\special{pn 8}%
\special{pa 2130 2710}%
\special{pa 2132 2742}%
\special{pa 2134 2774}%
\special{pa 2138 2806}%
\special{pa 2143 2837}%
\special{pa 2151 2868}%
\special{pa 2162 2899}%
\special{pa 2174 2928}%
\special{pa 2189 2957}%
\special{pa 2206 2984}%
\special{pa 2227 3009}%
\special{pa 2251 3031}%
\special{pa 2278 3048}%
\special{pa 2309 3061}%
\special{pa 2340 3069}%
\special{pa 2372 3070}%
\special{pa 2403 3065}%
\special{pa 2435 3061}%
\special{pa 2469 3060}%
\special{pa 2502 3058}%
\special{pa 2528 3043}%
\special{pa 2530 3040}%
\special{sp -0.045}%
%
\special{pn 8}%
\special{pa 2440 2440}%
\special{pa 2464 2419}%
\special{pa 2490 2399}%
\special{pa 2516 2382}%
\special{pa 2545 2368}%
\special{pa 2576 2359}%
\special{pa 2607 2355}%
\special{pa 2639 2357}%
\special{pa 2670 2363}%
\special{pa 2702 2373}%
\special{pa 2733 2385}%
\special{pa 2763 2400}%
\special{pa 2793 2416}%
\special{pa 2822 2433}%
\special{pa 2848 2451}%
\special{pa 2872 2472}%
\special{pa 2892 2496}%
\special{pa 2909 2524}%
\special{pa 2920 2555}%
\special{pa 2928 2587}%
\special{pa 2930 2619}%
\special{pa 2930 2651}%
\special{pa 2929 2683}%
\special{pa 2926 2715}%
\special{pa 2921 2747}%
\special{pa 2913 2778}%
\special{pa 2910 2809}%
\special{pa 2918 2840}%
\special{pa 2930 2870}%
\special{sp -0.045}%
\put(26.5000,-30.9000){\makebox(0,0)[lb]{$S^0\x B^2$}}%
%
\put(18.0000,-34.0000){\makebox(0,0)[lb]{}}%
%
\put(23.5000,-36.0000){\makebox(0,0)[lb]{}}%
\end{picture}%

\vskip-10mm
{\bf Figure \ref{high}.1:  The (1,2)-pass move changes 
a 2-dimensional closed oriented submanifolds $J\subset S^4$ into $K\subset S^4$.}
     
\bigbreak
\end{figure}

By Theorem \ref{r12}, 
$J$ and $K$ differ only by a single (1,2)-pass-move in a 4-ball $B^4$ 
embedded in $S^4$.   
By using Pontrjagin-Thom construction, we have the following Lemma \ref{mazu}. 
In Figure \ref{high}.1 we draw $J\cap B^4$ and $K\cap B^4$.

\begin{lem}\label{mazu} 
There is a Seifert hypersurface  
$V_J$ $($respectively, $V_K)$   
for $J$ $($respectively, $K)$ 
such that 
$V_J\cap B^4$ 
$($respectively, $V_K\cap B^4)$  
is \\
$($a 3-dimensional 1-handle$)\amalg($a 3-dimensional 2-handle$)$   
as in Definition \ref{lion}. 
\end{lem}

In Figure \ref{high}.2 we draw $V_J\cap B^4$ and $V_K\cap B^4$.

\begin{figure}
     \begin{center}

\input Figure1.tex
     \end{center}

\vskip-12mm
{\bf Figure \ref{high}.2:  
The (1,2)-pass move changes the Seifert hypersurfaces $V_J$ for $J$ 
into  $V_K$ for $K$.}

\bigbreak
\end{figure}

\bigbreak \noindent
{\bf Remark on Handle Notation.} 
Take a handle decomposition of $V_J$ (respectively, $V_K$)    
which consists of a single 3-dimensional 0-handle $h^0_J$  (respectively, $h^0_K$),    
the above 3-dimensional 1-handle, the above 3-dimensional 2-handle and other handles. 
Call this 3-dimensional 1-handle $h^1_J$ (respectively, $h^1_K$).  
Call this 3-dimensional 2-handle $h^2_J$ (respectively, $h^2_K$).  
We abbreviate by removing the subscript $_J$ (respectively, $_K$) when this is clear from the context.

\begin{lem}\label{persimmon} 
Under the assumption of Lemma \ref{mazu}, 
we can suppose that the attached part of the 3-dimensional 1-handle $h^1$ 
$($respectively, 2-handle $h^2)$ is embedded in the boundary of the 3-dimensional 0-handle $h^0$.  
\end{lem}

\noindent{\bf Note.}
Take a handle decomposition of a connected compact manifold. 
Let $i$ be a positive integer.
It does not hold in general that 
the attached part of each $i$-handle is embedded in the boundary of a 0-handle.

\bigbreak
\noindent{\bf Proof of Lemma \ref{persimmon}.}  
Take the tubular neighborhood $N(B^4)$ of $B^4$ in $S^4$. 
Note suppose that there exists an orientation preserving diffeomorphism map 
$\nu:N(B^4)\to B^4$ 
such that $\nu\vert_{V\cap N(B^4)}:V\cap N(B^4)\to V\cap B^4$ 
is an orientation preserving diffeomorphism map. 
Note $h^2$ in Lemma \ref{mazu}.  
Take a 4-dimensional 2-handle $H^2$ embedded in $($Int$N(B^4))-B^4$ and attach it to $V$ 
so that the attached part of $H^2$ is contained in $h^2.$   
See Figure \ref{high}.3. 
This surgery makes a new Seifert hypersurface, and call it $V$ again.

\begin{figure}
     \begin{center}
\input tubularx.tex
     \end{center}
\bigbreak
\end{figure}

Let $*=1,2$ and let $h^*\cap B^4$ be a new $h^*$. 
Take a  3-dimensional 0-handle $h^0$ in 
$(\partial(H^3))-$Int(the attached part) 
so that `the attached part of $h^2$' is embedded in $h^0$.



By using an isotopy and cancellation of handles in $V$,  we can assume that \newline 
`(the attached part of $h^1$) $\amalg$ (that of $h^2$)' 
is embedded in $h^0$ (respectively, $\partial B^4$) 
and that $h^1\amalg h^2$ is embedded in $B^4$. 

Thus we obtain a handle decomposition that satisfies the condition of Lemma \ref{persimmon}. \qed

\bigbreak
Recall the empty knot $[2]$, and the relation, Hopf link$=[2]\otimes[2]$, 
in Note \ref{kouiunomo}.




By the definition of knot products, we have the following. 

\begin{lem}  \label{soda}  
{\rm{\bf(This lemma follows from \cite[Lemma 3.4]{KauffmanNeumann} 
which is cited in \S\ref{amain}.)}} 
Let $p$ be a nonnegative integer. 
Let $E$ be a $p$-dimensional $($not necessarily connected or spherical$)$ 
closed oriented submanifold contained in $S^{p+2}$. 
The $(p+2)$-dimensional closed oriented submanifold 
$E\otimes[2]$ in the standard $(p+4)$-sphere 
is constructed as follows: 

Take the standard $(p+2)$-sphere $S^{p+2}$ in 
the standard $(p+4)$-sphere $S^{p+4}$ in the standard position. 
Take $E$ in $S^{p+2}$. 
There is a submanifold ${S'}^{p+2}$ of $S^{p+4}$ with the following properties: 
$S^{p+2}$ is ambient isotopic to ${S'}^{p+2}$, keeping $E$. 
%
%
%
%
$S^{p+2}$ intersects ${S'}^{p+2}$ transversely at $E$.  
Note that $S^{p+2}\cap {S'}^{p+2}$ in $S^{p+2}$ $($respectively, ${S'}^{p+2})$ is the submanifold $E$ of $S^{p+2}$ $($respectively, ${S'}^{p+2})$. 

Take the double branched cyclic covering space of $S^{p+4}$ along $S^{p+2}$. 

Then $S^{p+4}$ becomes the standard $(p+4)$-sphere again, and call it ${S'}^{p+4}$. 
${S'}^{p+2}$ becomes a closed oriented submanifold contained in ${S'}^{p+4}$. 
This submanifold is $E\otimes[2]\subset{S'}^{p+4}$.

$$
\begin{CD} 
E\otimes[2]&@>embedding>>&{S'}^{p+4} \\
@V{\text{branched along $E$}}VV& @VV{\text{branched along $S^{p+2}$}}V \\
{S'}^{p+2} &@>>embedding>& {S}^{p+4}
  \end{CD}
$$
\end{lem}

\noindent{\bf Proof of Lemma \ref{soda}.}  
Let $\nu$ be a positive integer. 
Define $\tau_\nu:D^2\to D^2$ as follows: 
Let $D^2=\{(a,b)|a^2+b^2\leqq1\}$. 
Let $\tau_\nu((0,0))=(0,0)$. 
Let $\theta\in\R$. Let $0<r\leqq1$. 
Let $\tau_\nu((r\cdot{\rm cos}\theta, r\cdot{\rm sin}\theta))=
(r\cdot{\rm cos}\hskip1mm (\nu\cdot\theta),r\cdot{\rm sin}\hskip1mm (\nu\cdot\theta))$.

Take the diagram in \cite[Lemma 3.4]{KauffmanNeumann}, which is cited in \S\ref{amain}.

$$
\begin{CD}
b_\nu(D^{p+3},F)&@>>>&D^{2}\\
@VVV& &@VV{\tau_\nu}V\\
D^{p+3}&@>{\gamma}>>&D^2
\end{CD}
$$
\\

Recall that $\gamma^{-1}((0,0))=F$ and that 
the submanifold $\partial F\subset \partial D^{p+3}$ is equal to 
the submanifold $E$ of the standard $(p+2)$-sphere. 

Recall that $b_\nu(D^{p+3},F)=\{(x,(a,b))\in D^{p+3}\x D^2\vert \gamma(x)=\tau_\nu((a,b))\}$, 
and that 
$\partial(b_\nu(D^{p+3},F))\subset \partial(D^{p+3}\x D^2)$ is 
the submanifold $E\otimes[\nu]\subset S^{p+4}$.

Let $\sigma_2: D^{p+3}\x D^2\to  D^{p+3}\x D^2$ be a map $({\rm id}, \tau_2)$. 
Note that the map $\sigma_2$ is 
a double branched cyclic covering whose branch set is $D^{p+3}\x\{(0,0)\}$. 
Note that 
$\partial(D^{p+3}\x\{(0,0)\})\\\subset \partial(D^{p+3}\x D^2)$
is 
the standard position of $(p+2)$-sphere in the standard $(p+4)$-sphere, 
and hence is the trivial $(p+2)$-knot, 
and let it be $S^{p+2}$. 
Note that 
Im$\sigma_2|_{\partial(D^{p+3}\x D^2)}$ is $\partial(D^{p+3}\x D^2)$, 
and that $\sigma_2|_{\partial(D^{p+3}\x D^2)}:{\partial(D^{p+3}\x D^2)}\to{\partial(D^{p+3}\x D^2)}$
 is a double branched cyclic cover whose branch set is $\partial(D^{p+3}\x\{(0,0)\})$ which is $S^{p+2}$. 

Consider 
$b_1(D^{p+3},F)\subset D^{p+3}\x D^2$ 
and 
$b_2(D^{p+3},F)\subset D^{p+3}\x D^2$. 
Note that \\
$\sigma_2^{-1}(b_1(D^{p+3},F))=b_2(D^{p+3},F)$  
and that 
$\sigma_2: b_2(D^{p+3},F)\to b_1(D^{p+3},F)$ is 
a double branched cyclic covering whose branch set is $F\x\{(0,0)\}$.

Since $\tau_1((a,b))=(a,b)$, we have 
$b_1(D^{p+3},F)=\{(x,(a,b))\in D^{p+3}\x D^2\vert \gamma(x)=(a,b)\}.$
Hence $(b_1(D^{p+3},F), D^{p+3}\x D^2)$ is a trivial pair of  (the $(p+3)$-ball, the $(p+5)$-ball) 
by a diffeomorphism of the $(p+5)$-ball. 
Hence $\partial(b_1(D^{p+3},F))\subset \partial(D^{p+3}\x D^2)$ is the trivial $p$-knot, 
and let it be ${S'}^{p+2}$.

Note the following: 
in $D^{p+3}\x D^2$, we have $(D^{p+3}\x\{(0,0)\})\cap b_1(D^{p+3},F)=F\x\{(0,0)\}.$   
$D^{p+3}\x\{(0,0)\}$ 
is ambient isotopic to $b_1(D^{p+3},F)$ in $D^{p+3}\x D^2$, keeping $F\x\{(0,0)\}.$   
By the definition of $\gamma$, $D^{p+3}\x\{(0,0)\}$ intersects $b_1(D^{p+3},F)$ transversely.

Therefore we have that, 
in $\partial(D^{p+3}\x D^2)$, we have 
$S^{p+2}\cap{S'}^{p+2}=\partial(F\x\{(0,0)\}).$   

Recall that 
the submanifold $\partial(F\x\{(0,0)\})\subset \partial(D^{p+3}\x\{(0,0)\})$ is equal to 
the submanifold $E\subset S^{p+2}$.

Hence we have the following: 
${S}^{p+2}$ is ambient isotopic to ${S'}^{p+2}$ in $S^{p+4}$, keeping \\$S^{p+2}\cap{S'}^{p+2}$.  
$S^{p+2}\cap{S'}^{p+2}$  in  
$S^{p+2}$ (respectively, ${S'}^{p+2}$)   
is the submanifold $E$ of $S^{p+2}$ (respectively, ${S'}^{p+2}$).

Therefore Im $\sigma_2|_{\partial(b_2(D^{p+3},F))}$ is $S'^{p+2}$,  
and that 
$\sigma_2|_{\partial(b_2(D^{p+3},F))}:{\partial(b_2(D^{p+3},F))}\to S'^{p+2}$ 
is a branched cyclic covering whose branch set is $E$. 

Recall that 
$\partial(b_2(D^{p+3},F))\subset \partial(D^{p+3}\x D^2)$ is 
the submanifold $E\otimes[2]\subset {S'}^{p+4}$.

Hence Lemma \ref{soda} holds. \qed\\

We have the following lemma.  
We use properties of branched cyclic covering spaces written in \cite[\S XII]{Kauffmanon}.

\begin{lem}\label{sodasui}
{\rm{\bf(This lemma follows from \cite[Lemma 3.4]{KauffmanNeumann} 
which is cited in \S\ref{amain},  
\cite[\S XII]{Kauffmanon}, and Lemma \ref{soda}.)}} 
Take the manifolds and the submanifolds of Lemma $\ref{soda}$. 

\smallbreak\noindent$(1)$ 
A Seifert hypersurface $W$ for $E\otimes[2]$ is constructed as follows 
$($Note $E\otimes[2]\subset W\newline\subset{S'}^{p+4})$. 
    
Embed a $(p+3)$-ball ${D'}^{p+3}$ in $S^{p+4}$ 
which is a Seifert hypersurface for ${S'}^{p+2}$ and 
which intersects $S^{p+2}$ transversely. 
When we take the double branched cyclic covering space of ${S'}^{p+2}$ 
defined in Lemma $\ref{soda}$, 
${D'}^{p+3}$ becomes a Seifert hypersurface $W$ for $E\otimes[2]$.

\smallbreak\noindent$(2)$ 
Furthermore $W$ has the following properties.  

The map $W\to {D'}^{p+3}$ 
is 
a double branched cyclic covering whose branch set is $V\\={D'}^{p+3}\cap {S}^{p+2}$.
Note that $\partial V=E=S^{p+2}\cap {S'}^{p+2}$ 
and that $V(\subset S^{p+2})$ is 
a Seifert hypersurface  for $E$. 

$V$ has the following property:  
Take a Seifert hypersurface $V' (\subset {S'}^{p+2})$  for $E$. 
%
$V'$ and $V$ are moved each other in ${D'}^{p+3}$ by an isotopy, 
keeping $\partial V'=\partial V=E$.

$$
\begin{CD} 
E\otimes[2]=\partial W&@>embedding>>&\hskip-9mm W&@>embedding>>&{S'}^{p+4} \\
@V{\text{branched along $E$}}VV
&\hskip10mm@V{\text{branched along $V$}}VV 
&  \hskip25mm@VV{\text{branched along $S^{p+2}$}}V \\
{S'}^{p+2}=\partial{D'}^{p+3}&@>>embedding>&\hskip-9mm{D'}^{p+3}&@>>embedding>& {S}^{p+4}
  \end{CD}
$$

Thus the manifold $W$ is constructed as follows: 
Take the tubular neighborhood $N(V')$ of $V'$ in ${S'}^{p+2}$. 
Take two copies of ${D'}^{p+3}$. 
Attach the two copies of ${D'}^{p+3}$ by identifying two $N(V')$ by the identity map.   
\end{lem}

\noindent{\bf Note.}  
We use $B^\sharp$ for the ball where the local-moves are carried out. 
So we use $D^\natural$ for the above ball.


\bigbreak 

By the assumption, 
$J(\subset S^4)$ and $K(\subset S^4)$ 
differ only in the interior of the 4-ball $B^4\\\subset S^4$.  
Furthermore recall the following. 
There is a Seifert hypersurface $V_J$ (respectively, $V_K$) for $J$ (respectively, $K$) 
with the following properties: 

\smallbreak\noindent$(1)$ 
$V_J$ and $V_K$ differ only in the interior of the 4-ball $B^4\subset S^4$.

\smallbreak\noindent$(2)$ 
 $V_Z\cap B^4$ ($Z=J,K$) is the disjoint union of 
a 1-handle $h^1_Z$ and a 2-handle $h^2_Z$. 

\smallbreak\noindent$(3)$ 
 $h^*_Z$ ($*=1,2$) is attached to $V_Z-$Int$B^4$ by the trivial framing. 
\\

By \cite[Lemma 3.4]{KauffmanNeumann}, \cite[\S XII]{Kauffmanon}, 
Lemmas \ref{soda} and \ref{sodasui}, we have the following. 

\begin{cla}\label{sui}
$(1)$ 
There is a Seifert hypersurface $V_{J\otimes[2]}$ $($respectively, $V_{K\otimes[2]})$ 
for $J\otimes[2]$ $($respectively, $K\otimes[2])\subset S^6$ 
such that 
$V_{J\otimes[2]}$ and $V_{K\otimes[2]}$ differ only in the 6-ball $B^6\subset S^6$, 
and such that  
$V_{Z\otimes[2]}\cap B^6$ 
is $h^2_{Z\otimes[2]}\amalg h^3_{Z\otimes[2]}$, 
where  $Z=J,K$, and 
$h^*_{Z\otimes[2]} (*=2,3)$ is a 5-dimensional $*$-handle attached to 
$V_{Z\otimes[2]}-${\rm Int}$B^6$ 
whose attached part is embedded in $V_{Z\otimes[2]}\cap B^6$.

\smallbreak\noindent$(2)$ 
$V_{Z\otimes[2]}$ is made from $V_Z$ as written in Lemma \ref{sodasui}. 
Suppose that \\
$V_Z=h^0_Z\cup...\cup h^1_{Z,1}...\cup h^1_{Z,a}\cup h^2_{Z,1}...\cup h^2_{Z,b}$ 
is a handle decomposition of $V_Z$, 
where $a$ and $b$ are positive integers, and $h^*_{Z,1}=h^*_Z  (*=1,2)$. 
Then 
$V_{Z\otimes[2]}$ has a handle decomposition \\
$V_{Z\otimes[2]}=
h^0_{Z\otimes[2]}\cup...\cup h^2_{Z\otimes[2],1}...\cup h^2_{Z\otimes[2],a}
\cup h^3_{Z\otimes[2],1}...\cup h^3_{Z\otimes[2],b}$ 
such that 
$h^{*+1}_{Z\otimes[2],\#}$ $(\#=1,...,a$ if $*=1$ and $\#=1,...,b$ if $*=2)$ is made from 
$h^*_{Z,\#}$ as written in  \cite[\S XII]{Kauffmanon}, 
 and such that $h^{*+1}_{Z\otimes[2],1}=h^{*+1}_{Z\otimes[2]}$.  
\end{cla}

Recall that, 
by Lemma \ref{mazu},  in $V_Z$, 
$h^1_Z$ and $h^2_Z$ is attached to $h^0_Z$ by the trivial framing. 
Hence we have the following.

\begin{cla}\label{kin}
Let $*=1,2$. 

\smallbreak\noindent$(1)$ 
We have that $h^*_Z\cup h^0_Z$ is $S^*\x($the $(3-*)$-ball$)$. 

\smallbreak\noindent$(2)$ 
Regard $h^*_Z\cup h^0_Z$ in $(1)$ as the trivial $(3-*)$-ball-bundle $\xi$ over $S^*$. 
Regard $S^4$ as $\partial D^5$ as written in Lemma \ref{sodasui}. 
Regard $S^*\x \{0\}$ as the boundary of an $(*+1)$-ball $\check B^{*+1}$ 
embedded in $\subset D^5$.  
Note that there is a framing $s$ on $\xi$ with the following properties: 
Restrict the normal bundle $h^*_Z\cup h^0_Z$ of $S^4$ to $S^*\x \{0\}$, and call the bundle $\varepsilon$. 
Make a framing on $\xi\oplus\varepsilon$ by using $s$ and the trivial framing on $\varepsilon$, 
and call the framing $t$.  
Take the normal bundle $\nu$ of $\check B^{*+1}$ of $D^5$. 
Note that $\nu|_{\partial \check B^{*+1}}$ is $\xi\oplus\varepsilon$. 

Then we have the following: The framing $t$ is extended over $\nu$.\\ 

\smallbreak\noindent$(3)$
$h^{*+1}_{Z\otimes[2]}$ is attached to $h^0_{Z\otimes[2]}$ by the trivial framing. 

Note that the union of two copies of $\check B^{*+1}$ is the core of 
$h^{*+1}_{Z\otimes[2]}\cup h^0_{Z\otimes[2]}$, where the two copies are made 
when we take the double branched cyclic cover in Lemma \ref{sodasui}. 
\end{cla}
\bigbreak

\begin{defn}\label{calH}
Let $n$ and $l$ be nonnegative integers. 
Let $P$ be an $n$-dimensional closed oriented submanifold of $S^{n+2}$. 
Let $X$ be a Seifert hypersurface for $P$. 
Let $B$ be the $(n+2)$-ball embedded in $S^{n+2}$.  
Assume that $X$ has a handle decomposition with the following properties: 
It consists of only one 0-handle $h^0$, 
and $l$-handles if $n=2l-1$
(respectively, $(l+1)$-handles and $(l+2)$-handles if $n=2l$).  

Suppose the following: $X\cap B$ is a disjoint union of two $l$-handles if $n=2l-1$ 
(respectively, an $(l+1)$-handle $h^{l+1}$ and an $(l+2)$-handle $h^{l+2}$ if $n=2l$).  
The attaching part of the two handles are embedded in $\partial h^0$. 
Note that the attaching parts are embedded in $\partial B$. 
These two handles are attached by the trivial framing. 

We say that 
$P$, $X$, the two handles, and $B$  
have a {\it property $\mathcal H$} 
if the following two bundles are equivalent. 

\smallbreak\noindent(1)
Suppose that each of the two handles is attached to $\partial h^0$ with the trivial framing. 
Let $h^k$ be the $k$-handle. 
Regard $h^k\cup h^0$ as 
$S^k\x($the $(n-k)$-ball$)$. 
Suppose that $S^k\x\{*\}$ 
is bounded by the $(k+1)$-ball $\hat B$ embedded in  $S^{n+2}$. 
The outer vector of $\hat B$ makes a bundle over  $S^k\x\{*\}$. 

\smallbreak\noindent(2) 
Restrict the normal bundle of $h^k\cup h^0$ in $S^{n+2}$ to $S^k\x\{*\}$. 
\end{defn}

\noindent{\bf Note.}  
Although `the two bundles (1) and (2) over $S^k\x\{*\}$ inDefinition \ref{calH}' 
may not be equivalent,  
the property $\mathcal H$ requests that the two bundles are equivalent. 
\bigbreak

By \cite[Lemma 6.1, its proof and section 6]{KauffmanNeumann}, we have the following.

\begin{cla}\label{dankai} 
$K\otimes[2]$, $Z\otimes[2]$, $h^2_{Z\otimes[2]}$, $h^3_{Z\otimes[2]}$ and $B^6$ 
in Claim $\ref{kin}$ 
have the property $\mathcal H$.
\end{cla}

By \cite[Lemma 3.4]{KauffmanNeumann}, \cite[\S XII]{Kauffmanon}, 
Lemmas \ref{soda}, \ref{sodasui}, and 
 \cite[Lemma 6.1, its proof and section 6]{KauffmanNeumann}, we have the following. 
We use a mathematical induction on $m$. 

\begin{cla}\label{dansui} 
Let $m$ be a nonnegative integer. 

\smallbreak\noindent 
$(1)$ 
There is a Seifert hypersurface $V_{J\otimes^m[2]}$ $($respectively, $V_{K\otimes^m[2]})$ 
for $J\otimes^m[2]$ $($respectively, $K\otimes^m[2])\subset S^{2m+4}$ 
such that 
$V_{J\otimes^m[2]}$ and $V_{K\otimes^m[2]}$ differ only 
in the $(2m+4)$-ball $B^{2m+4}\\\subset S^{2m+4}$, 
and such that  
$V_{Z\otimes^m[2]}\cap B^{2m+4}$ $(Z=J,K)$ 
is $h^{m+1}_{Z\otimes[2]}\amalg h^{m+2}_{Z\otimes^m[2]}$, 
and 
$h^{m+*}_{Z\otimes^m[2]} (*\\=1,2)$ is a $(2m+3)$-dimensional $(m+*)$-handle attached to 
$V_{Z\otimes^m[2]}-${\rm Int}$B^{2m+4}$ 
whose attached part is embedded in $V_{Z\otimes^m[2]}\cap B^{2m+4}$.

\smallbreak\noindent$(2)$ 
$V_{Z\otimes^{m+1}[2]}$ is made from $V_{Z\otimes^m[2]}$ as written in Lemma \ref{sodasui}. 
Suppose that \\
$V_{Z\otimes^m[2]}=
h^0_{Z\otimes^m[2]}\cup...\cup h^{m+1}_{Z\otimes^m[2],1}...\cup h^{m+1}_{Z\otimes^m[2],a}
\cup h^{m+2}_{Z\otimes^m[2],1}...\cup h^{m+2}_{Z\otimes^m[2],b}$ 
is a handle decomposition of $V_{Z\otimes^m[2]}$, 
where $a$ and $b$ are positive integers, and 
$h^{m+*}_{Z\otimes^m[2],1}=h^{m+*}_{Z\otimes^m[2]}  (*=1,2)$. 

Then 
$V_{Z\otimes^{m+1}[2]}$ has a handle decomposition \\
$V_{Z\otimes^{m+1}[2]}=
h^0_{Z\otimes^{m+1}[2]}\cup...\cup h^{m+2}_{Z\otimes^{m+1}[2],1}...\cup h^{m+2}_{Z\otimes^{m+1}[2],a}
\cup h^{m+3}_{Z\otimes^{m+1}[2],1}...\cup h^{m+3}_{Z\otimes^{m+1}[2],b}$ 
such that 
$h^{m+2}_{Z\otimes^{m+1}[2],\#}$ $(\#=1,...,a)$ is made from $h^{m+1}_{Z\otimes^m[2],\#}$ 
as written in  \cite[\S XII]{Kauffmanon},  
 such that 
$h^{m+3}_{Z\otimes^{m+1}[2],\#}$ $(\#=1,...,b)$ is made from $h^{m+2}_{Z\otimes^m[2],\#}$, 
and such that \\
$h^{m+1+*}_{Z\otimes^{m+1}[2],1}=h^{m+1+*}_{Z\otimes^{m+1}[2]}  (*=1,2)$. 


\smallbreak\noindent$(3)$
$Z\otimes^m[2]$, $V_{Z\otimes^m[2]}$, 
$h^{m+1}_{Z\otimes^m[2]}$, $h^{m+2}_{Z\otimes^m[2]}$, and $B^{2m+4}$ 
have the property $\mathcal H$.
\end{cla}





We have the following. 

\begin{cla}\label{chi}
Let $W^3$ be a submanifold $V_J-${\rm Int}$B^4$ of $S^4$. 
Note that the submanifold \\$V_K-${\rm Int}$B^4$ of $S^4$ is also $W^3$. 
Let $A^2$ be a closed oriented 2-dimensional submanifold $\partial W^3$ of $S^4$. 
Call $V_{Z\otimes[2]}-${\rm Int}$B^6$, $W^5$.  
Let $A^4$ be a closed oriented 4-dimensional submanifold $\partial W^5$ of $S^6$. 
We have that $A^4=A^2\otimes[2]$ and that $A^4$ is simply connected.  
\end{cla}






We have the following. 

\begin{cla}\label{moku}
Let $m$ be a positive integer.
Let $W^{2m+3}$ be a submanifold $V_{Z\otimes^m[2]}- B^{2m+4}$ of $S^{2m+4}$. 
Let $A^{2m+2}$ be a submanifold $\partial W^{2m+3}$ of $S^{2m+4}$.   
We have that 
$A^{2m+2}=A^2\otimes^m[2]$ 
and that 
$A^{2m+2}$ is $m$-connected. 
\end{cla}

Let $Q$ be a compact manifold with a handle decomposition. 
Let $h^q$ be a $q$-handle in the handles of the handle decomposition.  
If the core of the attached part of $h^q$ is null-homologous 
in the $(q-1)$-handle body of the handle decomposition,  
$h^q$ represents an element in $H_q(Q; \Z)$, and let  $[h^q]$ denote the element.  

\bigbreak 
Let $\xi$ and $\zeta$ be cycles in a compact oriented manifold $R$. 
Let $\xi\cdot\zeta$ in $R$ 
denote the intersection product of $\xi$ and $\zeta$ in $R$. 
We sometimes delete `in $R$' when this is clear from the context.

\bigbreak
By Lemma \ref{persimmon},    
$[h^1_Z]$ and $[h^2_Z]$ make sense. 
Recall that the attached part of  
$h^{m+\varepsilon}_{Z\otimes^m[2]} (\varepsilon\\=1,2)$    
is embedded in the $(2m+3)$-dimensional 0-handle in $V_{Z\otimes^m[2]}$.   
Hence 
$[h^{m+1}_{Z\otimes^m[2]}]$ and $[h^{m+2}_{Z\otimes^m[2]}]$ make sense. 
By the assumption of Main Theorem \ref{melon} 
$$
[h^1_{J}]\cdot[h^2_{J}] 
\text{ in $V^3_{J}$}
= 
[h^1_{K}]\cdot[h^2_{K}] 
\text{ in $V_{K}$}.  
$$

By \cite[Proposition 6.2 and \S6]{KauffmanNeumann}  
we have the following lemma. 

\begin{lem}\label{banana}   
Let $k$ be a positive integer, and $m=2k$.    
Let $Z=J,K$. 
We have
$$
[h^{m+1}_{Z\otimes^m[2]}]\cdot[h^{m+2}_{Z\otimes^m[2]}] 
\text{ in $V_{Z\otimes^m[2]}$}
=
[h^1_Z]\cdot[h^2_Z] 
\text{ in $V_Z$}. 
$$
\end{lem}

\noindent
{\bf Note.}  Recall $Z\otimes^{2k}[2]=Z\otimes^k\mathrm{Hopf}$.  
Note the difference between Lemma \ref{banana} and Lemma \ref{pear}. 

\bigbreak
By this lemma, we have the following: if $m=2k$ and   $k$ is a positive integer,   
$$
[h^{m+1}_{J\otimes^m[2]}]\cdot[h^{m+2}_{J\otimes^m[2]}] 
\text{ in $V_{J\otimes^m[2]}$}
=
[h^{m+1}_{K\otimes^m[2]}]\cdot[h^{m+2}_{K\otimes^m[2]}] 
\text{ in $V_{K\otimes^m[2]}$}. 
$$

Therefore we have  the following. 
\begin{cla}\label{do}
Let $\check {S}_Z^m$ $($respectively, $\check {S}_Z^{m+1})$ be
 the standard $m$-$($respectively, $(m+1))$-sphere
which is the core of the attached part of the handle $h^{m+1}_{Z\otimes^m[2]}$ $($respectively, $h^{m+2}_{Z\otimes^m[2]})$.   
Then there is a diffeomorphism map 

\noindent
$f:V_{J\otimes^m[2]}-${\rm Int}$B^{2m+4}\to V_{K\otimes^m[2]}-${\rm Int}$B^{2m+4}$ 

\noindent
such that 
$f|_{\partial(V_{J\otimes^m[2]}-{\rm Int}B^{2m+4})}$ 
 carries 
the homology class of $\check {S}_J^{m}$ $($respectively, $\check {S}_J^{m+1})$     
to that of  $\check {S}_K^{m}$ $($respectively, $\check {S}_K^{m+1}).$
\end{cla}

Note that the submanifold  $\partial(V_{J\otimes^m[2]}-{\rm Int}B^{2m+4})$ of $S^{2m+4}$ is $A^{2m+2}$, which is defined in Claim \ref{moku}.
\\

By \cite{Haefligerunknot,  Whitney, Whitneytrick}, 
we have the following. 

\begin{cla}\label{ten}
Let $m\geqq4$. 
 The embedding type of 
$\check {S}^{m}$ $($respectively, $\check {S}^{m+1})$ 
in $A^{2m+2}$ is unique 
if the homology class 
$[\check {S}^{m}]$ $($respectively, $[\check {S}^{m+1}])$ 
in $A^{2m+2}$ is fixed. 
\end{cla}

Therefore we have the following. 

\begin{cla}\label{mei}
There is a diffeomorphism map 
from 
$V_{J\otimes^m[2]}-${\rm Int}$B^{2m+4}$ 
to 
$V_{K\otimes^m[2]}-${\rm Int}$B^{2m+4}$ 
such that 
the image of 
the attached part of $h^{m+*}_{J\otimes^m[2]}$
is 
that of $h^{m+*}_{K\otimes^m[2]}$ $(*=1,2)$. 
\end{cla}

Recall that $h^{m+*}_{Z\otimes^m[2]}$ ($*=1,2. Z=J,K$)
is attached by the trivial framing to the only one 0-handle in 
$V_{Z\otimes^m[2]}-${\rm Int}$B^{2m+4}$. 
Therefore 
the diffeomorphism in Claim \ref{mei} extends to one 
from $V_{J\otimes^m[2]}$ to $V_{K\otimes^m[2]}$.  
 

\bigbreak
By \cite[Proposition 6.2 and \S6]{KauffmanNeumann}  
we have the following lemma. 

\begin{lem} \label{pineapple}
Let $m=2k$ and  $k$ a positive integer. 
The $\Z$-Seifert pairing of $[h^{m+1}_{Z\otimes^m[2]}]$ and $[h^{m+2}_{Z\otimes^m[2]}]$  
associated with the Seifert hypersurface  $V_{Z\otimes^m[2]}$  is equal to \newline
$(-1)^k\x$$($that of $[h^1_Z]$ and $[h^2_Z] $ associated with $V_Z)$.  
\end{lem}

Note that the handles, $h^{m+1}_{Z\otimes^m[2]}$ and $h^{m+2}_{Z\otimes^m[2]}$, 
embedded in $B^{2m+4}_{Z\otimes^{2m}}$ 
are attached to $\partial B^{2m+4}_{Z\otimes^{2m}}$. 
Note that  the embedding type of 
the disjoint union of the attached part of $h^{m+1}_{Z\otimes^m[2]}$ and 
that of $h^{m+2}_{Z\otimes^m[2]}$ in  $\partial B^{2m+4}_{Z\otimes^{2m}}$ 
is only one. 
\\

Recall that 
$h^{m+1}_{Z\otimes^m[2]}$ and $h^{m+2}_{Z\otimes^m[2]}$ are embedded in  $B^{2m+4}_{Z\otimes^{2m}}$. 
Let $\tilde{B}^{m+1}$ be the core of $h^{m+1}_{Z\otimes^m[2]}$.  
Note that we regard $[\tilde{B}^{m+1}, \partial(\tilde{B}^{m+1})]$ 
as an element in \\
$H_{m+1}(
\overline{B^{2m+4}_{Z\otimes^{2m}}-h^{m+2}_{Z\otimes^m[2]}}, 
\partial(
\overline{B^{2m+4}_{Z\otimes^{2m}}-h^{m+2}_{Z\otimes^m[2]}}
);
\Z)$. 
Here, we use the following notation: 
Let $E$ be a sub-topological space of a topological space $F$. 
Let $\overline{E}$ denote the closure of $E$ in $F$. 

By \cite[Lemma 3.4]{KauffmanNeumann}, \cite[\S XII]{Kauffmanon}, Lemmas \ref{soda} and \ref{sodasui},  
there is a one-to-one correspondence between the following two sets (1) and (2) 
if $m$ is an integer$\geq2$. 

\smallbreak\noindent (1) 
The set of the embedding type of 
the disjoint union of $h^{m+1}_{Z\otimes^m[2]}$ and $h^{m+2}_{Z\otimes^m[2]}$ in  $B^{2m+4}_{Z\otimes^{2m}}$,  keeping their attached part. 

\smallbreak\noindent (2) 
The set 
$H_{m+1}(
\overline{B^{2m+4}_{Z\otimes^{2m}}-h^{m+2}_{Z\otimes^m[2]}}, 
\partial(
\overline{B^{2m+4}_{Z\otimes^{2m}}-h^{m+2}_{Z\otimes^m[2]}}
);
\Z)$. 
\\

Furthermore there is a one-to-one correspondence between the above (2) and the following (3).

\smallbreak\noindent (3) 
The set of the $\Z$-Seifert pairing of $[h^{m+1}_{Z\otimes^m[2]}]$ and $[h^{m+2}_{Z\otimes^m[2]}]$  
associated with the Seifert hypersurface  $V_{Z\otimes^m[2]}$.  
\\

Therefore 
$V_{J\otimes^m[2]}$ and $V_{K\otimes^m[2]}$ differ by a single $(m+1,m+2)$-pass-move 
if $m=2k$, where $k$ is any integer$\geq2$.

Hence $J\otimes^k{\rm Hopf}$ is obtained from 
$K\otimes^k{\rm Hopf}$  by a single 
$(2k+1,2k+2)$-pass-move  
where $k$ is any integer$\geq2$.

This completes the proof of Main Theorem \ref{melon}. 
\qed

\bigbreak 
\noindent 
{\bf Note.}  
 \cite[Proposition 6.2 and \S6]{KauffmanNeumann},  
\cite[Proposition 5.4 and section 5]{KauffmanOgasa}, and 
\cite[section 7]{Levinecob}  imply 
 Lemma \ref{pear} and \ref{bean}.

\begin{lem} \label{pear}  
Let $k$ be a nonnegative integer. 
If $m=2k+1$, \newline
$[h^{m+1}_{Z\otimes^m[2]}]\cdot[h^{m+2}_{Z\otimes^m[2]}] 
\text{ in $V_{Z\otimes^m[2]}$}
\neq 
[h^1_Z]\cdot[h^2_Z] 
\text{ in $V_Z$}
 \quad (Z=J,K)$\quad
in general. 
\end{lem}

Note the difference between Lemma \ref{banana} and Lemma \ref{pear}.

\begin{lem} \label{bean}
Let $k$ be a positive integer.
$J\otimes^{2k+1}[2]$ and 
$K\otimes^{2k+1}[2]$ are not diffeomorphic or homeomorphic in general 
\end{lem}

By this lemma we have the following:  
Let $k$ be a positive integer. 
$J\otimes^{2k+1}[2]$ and  
$K\otimes^{2k+1}[2]$ are 
 not $(2k+2, 2k+3)$-pass-move equivalent in general.



\bigbreak
\section{
Proof of Main Theorem \ref{sugar}}\label{hightsuzuki}
\noindent   
The idea of the proof is similar  to that in \S\ref{high}: 
Take a Seifert hypersurface $V_J$ (respectively $V_K$) for $J$ (respectively $K$) 
such that $V_J$ and $V_K$ differ only in the 3-ball 
where the one pass-move is carried out. 
We make a Seifert hypersurface 
$V_{J\otimes^k{\rm Hopf}}$ (respectively $V_{K\otimes^k{\rm Hopf}}$) 
for $J\otimes^k{\rm Hopf}$ (respectively $K\otimes^k{\rm Hopf}$)  
from $V_J$ (respectively $V_K$) 
by using Lemmas \ref{soda} and \ref{sodasui}
such that 
$V_{J\otimes^k{\rm Hopf}}$ and $V_{K\otimes^k{\rm Hopf}}$ 
differ only by one  $(2k+1,2k+1)$-pass-move 
in a $(4k+3)$-ball of $S^{4k+3}$.  
Their difference implies that 
$J\otimes^k{\rm Hopf}$ and $K\otimes^k{\rm Hopf}$  
differ by one $(2k+1,2k+1)$-pass-move. 

Let $m$ be a nonnegative integer. 
We point out that  
there are not 
a Seifert hypersurface $V_{J\otimes^{2m+1}[2]}$ for $J\otimes^{2m+1}[2]$ 
and 
a Seifert hypersurface $V_{K\otimes^{2m+1}[2]}$ for $K\otimes^{2m+1}[2]$ 
which are diffeomorphic 
 in general  
although  
$V_{J\otimes^k{\rm Hopf}}=V_{J\otimes^{2m}[2]}$ 
and 
$V_{K\otimes^k{\rm Hopf}}=V_{K\otimes^{2m}[2]}$ 
have the above property.

\bigbreak
We prove Theorem \ref{B} by using Main Theorem \ref{sugar} before proving  
Main Theorem \ref{sugar}. 

\smallbreak
\noindent
{\bf Proof of Theorem \ref{B}.} 
We first prove Theorem \ref{B}.(1). 
Let $K_{(a,b)}$ denote the $(a,b)$ torus link. 
By Main Theorem \ref{sugar} 
$K_{(a,b)}\otimes^k\mathrm{Hopf}$ is 
$(2k+1, 2k+1)$-pass-move equivalent to 
$K_{(a', b')}\otimes^k\mathrm{Hopf}$. 
Hence 
$K_{(a,b)}\otimes^k\mathrm{Hopf}$ is 
diffeomorphic to 
$K_{(a', b')}\otimes^k\mathrm{Hopf}$. 
By (\cite[the sixth line from the bottom of page 1106]{Kauffman}, 
$K_{(a,b)}\otimes^k\mathrm{Hopf}$ (respectively, $K_{(a',b')}\otimes^k\mathrm{Hopf}$)  
is diffeomorphic to 
$\Sigma(a,b,\underbrace{2,...,2}_{2k})$
(respectively, $\Sigma(a',b',\underbrace{2,...,2}_{2k})$). 
Hence $\Sigma(a,b,\underbrace{2,...,2}_{2k})$ 
is diffeomorphic to $\Sigma(a',b',\underbrace{2,...,2}_{2k})$.

We prove Theorem \ref{B}.(3).  
By the assumption, $J$ and $K$ are pass-move equivalent. 
Therefore, by Main Theorem \ref{sugar},  
$J\otimes^k{\rm Hopf}$ is pass-move equivalent  to $K\otimes^k{\rm Hopf}$.   
Hence  $J\otimes^k{\rm Hopf}$ is diffeomorphic to $K\otimes^k{\rm Hopf}$.   
Note that 
 Theorem \ref{B}.(1) follows from  Theorem \ref{B}.(3).  

We prove Theorem \ref{B}.(2).  
By \cite[Exercise 8 in p.177-178]{Rolfsen} and  \cite[Theorem 10.4 in p.261]{Kauffmanon}, 
there are the torus knots, $K_{(a,b)}$ and $K_{(a',b')}$, such that 
the Arf invariant of $K_{(a,b)}$ is different from that of $K_{(a',b')}$. 
Hence $K_{(a,b)}$ is not pass-move equivalent to $K_{(a',b')}$ by 
\cite[Proposition 5.6 in page 77 and Corollary in page 260]{Kauffmanon}.  

By (\cite[the sixth line from the bottom of page 1106]{Kauffman}, 
$K_{(a,b)}\otimes$Hopf (respectively, $K_{(a',b')}\otimes$Hopf) is diffeomorphic to 
$\Sigma(a, b, 2, 2)$ (respectively, $\Sigma(a', b', 2, 2)$).  
By Theorem \ref{LE}, 
$K_{(a,b)}\otimes$Hopf (respectively, $K_{(a',b')}\otimes$Hopf) is PL homeomorphic to the standard 5-sphere. 
By   \cite{KervaireMilnor},  we have the following. 
if a smooth 5-manifold $M$ is orientation preserving PL homeomorphic to the  standard sphere, 
then $M$ is orientation preserving diffeomorphic to the  standard sphere. 
Hence $K_{(a,b)}\otimes$Hopf (respectively, $K_{(a',b')}\otimes$Hopf) is diffeomorphic to the standard 5-sphere. 
Hence $\Sigma(a,b,2,2)$ is diffeomorphic to $\Sigma(a',b',2,2)$. \qed

\bigbreak
\noindent{\bf Note.} 
The number of the connected components of a torus knot is one. 
The number of the connected components of a torus link is 
no less than one.

\bigbreak
\noindent {\bf Proof of Main Theorem \ref{sugar}.}  
We need the following lemmas. 

Let $V_J$ (respectively, $V_K$) be a Seifert hypersurface for $J$ (respectively, $K$). 
We can suppose that 
$V_J$ and $V_K$ differ only in the 3-ball $B^3$ as shown in Figure \ref{hightsuzuki}.1. 
We have the following by the definition of the pass-move. 

\begin{figure}
     \begin{center}

\input XFigure1.tex
\end{center}
\vskip3mm
\end{figure}

\begin{lem}\label{Xmazu}
There is a handle decomposition of $V_J$ $($respectively, $V_K)$ with the following properties:  

\smallbreak \noindent$(1)$   
$V_J\cap B^3$ $($respectively, $V_K\cap B^3)$  
is a disjoint union of two 2-dimensional 1-handles 
as shown in Figure \ref{hightsuzuki}.1.  
Let $Z=J,K$. 
Call one of these 2-dimensional 1-handles $h_{\alpha, Z}$ and the other $h_{\beta, Z}$. 
The handles $h_{\alpha, Z}$ and $h_{\beta, Z}$   are embedded 
in $B^3$ as shown in Figure \ref{hightsuzuki}.1.

\smallbreak \noindent$(2)$ 
There is a single 2-dimensional 0-handle $h^0_J$ $($respectively, $h^0_K)$ such that 
the attached part of $h_{\alpha, Z}$ $($respectively,  $h_{\beta, Z})$  
is embedded in $\partial(h^0_J)$ $($respectively, $\partial(h^0_K))$, and furthermore     
is embedded in $\partial B^3$ as shown in Figure \ref{hightsuzuki}.1.  
\end{lem}

By attaching 3-dimensional 1-handles embedded in $S^3$ as in Lemma \ref{persimmon}, 
we have the following.  

\begin{lem}\label{zerofr}
Let $\zeta=\alpha, \beta$. 
Hence a Seifert pairing of $[h_{\zeta, Z}]$ and itself is zero.
\end{lem}

By the assumption, 
$J(\subset S^3)$ and $K(\subset S^3)$ differ only in the interior of 
the 3-ball $B^3\subset S^3$. 
%
%
By \cite[Lemma 3.4]{KauffmanNeumann}, \cite[\S XII]{Kauffmanon}, 
Lemmas \ref{soda} and \ref{sodasui}, we have the following. 

\begin{cla}\label{Xsui}
$(1)$ 
There is a Seifert hypersurface $V_{J\otimes[2]}$ $($respectively, $V_{K\otimes[2]})$ 
for $J\otimes[2]$ $($respectively, $K\otimes[2])\subset S^5$ 
such that $V_{J\otimes[2]}$ and $V_{K\otimes[2]}$ differ only in the 5-ball $B^5\subset S^5$, 
and such that  $V_{Z\otimes[2]}\cap B^5$ is 
$h^2_{\alpha, Z\otimes[2]}\amalg h^2_{\beta, Z\otimes[2]}$, 
where  $Z=J,K$, and 
$h^2_{\zeta, Z\otimes[2]}$ $(\zeta=\alpha, \beta)$ is a 4-dimensional 2-handle attached to $V_{Z\otimes[2]}-${\rm Int}$B^5$ whose attached part is embedded in $V_{Z\otimes[2]}\cap B^5$.

\smallbreak\noindent$(2)$ 
$V_{Z\otimes[2]}$ is made from $V_Z$ as written in Lemma \ref{sodasui}. 
Suppose that \\
$V_Z=h^0_Z\cup...\cup h^1_{Z,1}...\cup h^1_{Z,a}$  
is a handle decomposition of $V_Z$, where $a$ is a positive integer, 
$h^1_{Z,1}=h^1_{\alpha,Z}$, and $h^1_{Z,2}=h^1_{\beta,Z}$. 
Then 
$V_{Z\otimes[2]}$ has a handle decomposition \\
$V_{Z\otimes[2]}=
h^0_{Z\otimes[2]}\cup...\cup h^2_{Z\otimes[2],1}...\cup h^2_{Z\otimes[2],a}$ 
such that 
$h^2_{Z\otimes[2],\nu}$ is made from $h^1_{Z,\nu}$ $(\nu\in\{1,...,a\})$ as written in  \cite[\S XII]{Kauffmanon}, 
and such that 
$h^2_{Z\otimes[2],1}=h^2_{\alpha,Z\otimes[2]}$ and $h^2_{Z\otimes[2],2}=h^2_{\beta,Z\otimes[2]}$. 
\end{cla}






By Lemmas \ref{Xmazu} and \ref{zerofr},  in $V_Z$, 
$h^1_{\alpha,Z}$ and $h^1_{\beta, Z}$ is attached to $h^0_Z$ by the trivial framing. 
By \cite[Lemma 6.1, its proof and section 6]{KauffmanNeumann}, 
we have the following.

\begin{cla}\label{Xdankai} 
$K\otimes[2]$, $Z\otimes[2]$, 
$h^2_{\alpha,Z\otimes[2]}$, $h^2_{\beta, Z\otimes[2]}$ and $B^5$ 
in Claim $\ref{Xsui}$ have the property $\mathcal H$.
\end{cla}

By \cite[Lemma 3.4]{KauffmanNeumann}, \cite[\S XII]{Kauffmanon}, 
Lemmas \ref{soda},  \ref{sodasui},   and   
\cite[Lemma 6.1, its proof and section 6]{KauffmanNeumann}
we have the following. 
We use a mathematical induction on $m$.

\begin{cla}\label{Xdansui} 
Let $m$ be a nonnegative integer. 

\smallbreak\noindent 
$(1)$ 
There is a Seifert hypersurface $V_{J\otimes^m[2]}$ $($respectively, $V_{K\otimes^m[2]})$ 
for $J\otimes^m[2]$ $($respectively, $K\otimes^m[2])\subset S^{2m+3}$ 
such that 
$V_{J\otimes^m[2]}$ and $V_{K\otimes^m[2]}$ differ only 
in the $(2m+3)$-ball $B^{2m+3}\\\subset S^{2m+3}$, 
and such that  
$V_{Z\otimes^m[2]}\cap B^{2m+3}$ $(Z=J,K)$ 
is $h^{m+1}_{\alpha,Z\otimes[2]}\amalg h^{m+1}_{\beta,Z\otimes^m[2]}$, 
where 
$h^{m+1}_{\gamma,Z\otimes^m[2]} (\gamma=\alpha,\beta)$ is a $(2m+2)$-dimensional $(m+1)$-handle 
attached to 
$V_{Z\otimes^m[2]}-${\rm Int}$B^{2m+3}$ 
whose attached part is embedded in $V_{Z\otimes^m[2]}\cap B^{2m+3}$.

\smallbreak\noindent$(2)$ 
$V_{Z\otimes^{m+1}[2]}$ is made from $V_{Z\otimes^m[2]}$ as written in Lemma \ref{sodasui}. 
Suppose that \\
$V_{Z\otimes^m[2]}=
h^0_{Z\otimes^m[2]}\cup...\cup h^{m+1}_{Z\otimes^m[2],1}...\cup h^{m+1}_{Z\otimes^m[2],a}$ 
is a handle decomposition of $V_{Z\otimes^m[2]}$, 
where $a$ is a positive integer, and 
$h^{m+1}_{Z\otimes^m[2],1}=h^{m+1}_{\alpha, Z\otimes^m[2]}$ 
and 
$h^{m+1}_{Z\otimes^m[2],2}=h^{m+1}_{\beta, Z\otimes^m[2]}$. 

Then 
$V_{Z\otimes^{m+1}[2]}$ has a handle decomposition \\
$V_{Z\otimes^{m+1}[2]}=
h^0_{Z\otimes^{m+1}[2]}\cup...\cup h^{m+2}_{Z\otimes^{m+1}[2],1}...\cup h^{m+2}_{Z\otimes^{m+1}[2],a}$ 
such that 
$h^{m+2}_{Z\otimes^{m+1}[2],\#}$ $(\#=1,...,a)$ is made from $h^{m+1}_{Z\otimes^m[2],\#}$ 
as written in  \cite[\S XII]{Kauffmanon},  
such that 
$h^{m+2}_{Z\otimes^{m+1}[2],1}=h^{m+2}_{\alpha, Z\otimes^{m+1}[2]}$,  
and such that 
$h^{m+2}_{Z\otimes^{m+1}[2],2}=h^{m+2}_{\beta, Z\otimes^{m+1}[2]}$. 


\smallbreak\noindent$(3)$  
$Z\otimes^m[2]$, $V_{Z\otimes^m[2]}$, 
$h^{m+1}_{\alpha, Z\otimes^m[2]}$, $h^{m+1}_{\beta, Z\otimes^m[2]}$, and $B^{2m+3}$
have the property $\mathcal H$.
\end{cla}

We have the following. 

\begin{cla}\label{Xchi}
Let $W^2$ be a submanifold $V_J-${\rm Int}$B^3$ of $S^3$. 
Note that the submanifold \\$V_K-${\rm Int}$B^3$ of $S^3$ is also $W^2$. 
Let $A^1$ be a closed oriented 1-dimensional submanifold $\partial W^2$ of $S^3$. 
Call $V_{Z\otimes[2]}-${\rm Int}$B^5$, $W^4$.  
Let $A^3$ be a closed oriented 3-dimensional submanifold $\partial W^4$ of $S^5$. 
We have that $A^3=A^1\otimes[2]$.
\end{cla}






We have the following. 

\begin{cla}\label{Xmoku}
Let $m$ be a positive integer.
Let $W^{2m+2}$ be a submanifold $V_{Z\otimes^m[2]}-${\rm Int}$B^{2m+3}$ of $S^{2m+3}$. 
Let $A^{2m+1}$ be a submanifold $\partial W^{2m+2}$ of $S^{2m+3}$.   
We have that $A^{2m+1}=A^1\otimes^m[2]$ 
and that $A^{2m+1}$ is $(m-1)$-connected. 
\end{cla}

By Lemma \ref{Xmazu},    
$[h^1_{\alpha,Z}]$ and $[h^1_{\beta,Z}]$ make sense. 
Recall that the attached part of $h^{m+1}_{\gamma,Z\otimes^m[2]} (\gamma\\=\alpha,\beta)$    
is embedded in the $(2m+2)$-dimensional 0-handle in $V_{Z\otimes^m[2]}$.   
Hence 
$[h^{m+1}_{\gamma,Z\otimes^m[2]}]$  makes sense. 

By the assumption of Main Theorem \ref{sugar}  
$$
[h^1_{\alpha,J}]\cdot[h^1_{\beta,J}] \text{ in $V^2_{J}$}
= 
[h^1_{\alpha,K}]\cdot[h^1_{\beta,K}] \text{ in $V^2_{K}$}. 
$$

By \cite[Proposition 6.2 and \S6]{KauffmanNeumann}  
we have the following lemma. 

\begin{lem}\label{Xbanana}   
Let $k$ be a positive integer, and $m=2k$.    
Let $Z=J,K$. 
We have
$$
[h^{m+1}_{\alpha, Z\otimes^m[2]}]\cdot[h^{m+1}_{\beta,Z\otimes^m[2]}] \text{ in $V_{Z\otimes^m[2]}$}
=
[h^1_{\alpha,Z}]\cdot[h^1_{\beta,Z}] \text{ in $V_Z$}. 
$$
\end{lem}

\noindent
{\bf Note.}  Recall $Z\otimes^{2k}[2]=Z\otimes^k\mathrm{Hopf}$.  
Note the difference between Lemma \ref{banana} and Lemma \ref{Xpear}. 

\bigbreak
By this lemma, we have the following: if $m=2k$ and   $k$ is a positive integer,   
$$
[h^{m+1}_{\alpha, J\otimes^m[2]}\cdot[h^{m+1}_{\beta, J\otimes^m[2]}] \text{ in $V_{J\otimes^m[2]}$}
=
[h^{m+1}_{\alpha, K\otimes^m[2]}\cdot[h^{m+1}_{\beta, K\otimes^m[2]}] \text{ in $V_{K\otimes^m[2]}$}. 
$$

Therefore we have  the following. 

\begin{cla}\label{Xdo} 
Let $Z=J,K.$
Let $\check {S}_{\gamma,Z}^m$ be the standard $m$-sphere 
such that 
the core of the handle $h^{m+1}_{\gamma, Z\otimes^m[2]}$ 
$(\gamma=\alpha, \beta)$.   

Then there is a diffeomorphism map 

\noindent
$f:V_{J\otimes^m[2]}-${\rm Int}$B^{2m+3}\to V_{K\otimes^m[2]}-${\rm Int}$B^{2m+3}$ 

\noindent
such that 
$f|_{\partial(V_{J\otimes^m[2]}-{\rm Int}B^{2m+3})}$ 
 carries 
the homology class of $\check {S}_{\gamma,J}^{m}$ 
to that of  $\check {S}_{\gamma,K}^{m}$. 
\end{cla}

Note that the submanifold  $\partial(V_{J\otimes^m[2]}-{\rm Int}B^{2m+3})$ of $S^{2m+3}$ is $A^{2m+1}$, which is defined in Claim \ref{Xmoku}.
\\

By \cite{Haefligerunknot,  Whitney, Whitneytrick}, we have the following. 

\begin{cla}\label{Xten}
Let $m\geqq1$.  the embedding type of $\check {S}_{\gamma,Z}^{m}$ in $A^{2m+1}$ 
is unique if the homology class $[\check {S}_{\gamma,Z}^{m}]$ of $A^{2m+1}$ is fixed. 
\end{cla}

Therefore we have the following. 

\begin{cla}\label{Xmei}
There is a diffeomorphism map 
from 
$V_{J\otimes^m[2]}-${\rm Int}$B^{2m+3}$ 
to 
$V_{K\otimes^m[2]}-${\rm Int}$B^{2m+3}$ 
such that 
the image of the attached part of $h^{m+1}_{\gamma, J\otimes^m[2]}$
is that of $h^{m+1}_{\gamma,K\otimes^m[2]}$ $(\gamma=\alpha,\beta)$. 
\end{cla}

Recall that $h^{m+1}_{\gamma, Z\otimes^m[2]}$ ($\gamma=\alpha, \beta. Z=J,K$)
is attached by the trivial framing to the only one 0-handle in 
$V_{Z\otimes^m[2]}-$Int$B^{2m+3}$. 
Therefore the diffeomorphism in Claim \ref{Xmei} extends to one 
from $V_{J\otimes^m[2]}$ to $V_{K\otimes^m[2]}$.  
 

\bigbreak
By \cite[Proposition 6.2 and \S6]{KauffmanNeumann}  
we have the following lemma. 

\begin{lem} \label{pineapple}
Let $m=2k$ and  $k$ a positive integer. 
Let $Z=J,K$. 
The $\Z$-Seifert pairing of 
$[h^{m+1}_{\alpha, Z\otimes^m[2]}]$ and $[h^{m+1}_{\beta, Z\otimes^m[2]}]$  
associated with the Seifert hypersurface  $V_{Z\otimes^m[2]}$  is equal to \newline
$(-1)^k\x$$($that of $[h^1_{\alpha, Z}]$ and $[h^1_{\beta, Z}] $ associated with $V_Z)$.  
\end{lem}

Note that the handles, 
$h^{m+1}_{\alpha, Z\otimes^m[2]}$ and $h^{m+1}_{\beta, Z\otimes^m[2]}$, embedded in 
$B^{2m+3}_{Z\otimes^{2m}}$ 
are attached to $\partial B^{2m+3}_{Z\otimes^{2m}}$. 
Note that  the embedding type of 
the disjoint union of the attached part of 
$h^{m+1}_{\alpha, Z\otimes^m[2]}$ and 
that of 
$h^{m+1}_{\beta, Z\otimes^m[2]}$ in  $\partial B^{2m+3}_{Z\otimes^{2m}}$ 
is only one. 
\\

Recall that 
$h^{m+1}_{\alpha, Z\otimes^m[2]}$ and $h^{m+1}_{\beta, Z\otimes^m[2]}$ are embedded in  $B^{2m+3}_{Z\otimes^{2m}}$. 
Let $\gamma=\alpha, \beta$.  
Let $\tilde{B}^{m+1}_\gamma$ be the core of $h^{m+1}_{\gamma, Z\otimes^m[2]}$.  
Note that we regard $[\tilde{B}^{m+1}_\alpha, \partial(\tilde{B}^{m+1}_\alpha)]$ 
as an element in \\
$H_{m+1}(
\overline{B^{2m+3}_{Z\otimes^{2m}}-h^{m+1}_{\beta, Z\otimes^m[2]}}, 
\partial(
\overline{B^{2m+3}_{Z\otimes^{2m}}-h^{m+1}_{\beta, Z\otimes^m[2]}}
);
\Z)$. 

By \cite[Lemma 3.4]{KauffmanNeumann}, \cite[\S XII]{Kauffmanon}, Lemmas \ref{soda} and \ref{sodasui},  
there is a one-to-one correspondence between the following two sets (1) and (2). 

\smallbreak\noindent (1) 
The set of the embedding type of 
the disjoint union of $h^{m+1}_{\alpha, Z\otimes^m[2]}$ and $h^{m+1}_{\beta, Z\otimes^m[2]}$ in  $B^{2m+3}_{Z\otimes^{2m}}$,  keeping their attached part. 

\smallbreak\noindent (2) 
The set 
$H_{m+1}(
\overline{B^{2m+3}_{Z\otimes^{2m}}-h^{m+1}_{\beta, Z\otimes^m[2]}}, 
\partial(
\overline{B^{2m+3}_{Z\otimes^{2m}}-h^{m+1}_{\beta, Z\otimes^m[2]}}
);
\Z)$. 
\\

Furthermore there is a one-to-one correspondence between the above (2) and the following (3).

\smallbreak\noindent (3) 
The set of the $\Z$-Seifert pairing of 
$[h^{m+1}_{\alpha, Z\otimes^m[2]}]$ and $[h^{m+1}_{\beta, Z\otimes^m[2]}]$  
associated with the Seifert hypersurface  $V_{Z\otimes^m[2]}$.  
\\

Therefore 
$V_{J\otimes^m[2]}$ and $V_{K\otimes^m[2]}$ differ by a single $(m+1,m+1)$-pass-move 
if $m=2k$.

Hence $J\otimes^k{\rm Hopf}$ is obtained from 
$K\otimes^k{\rm Hopf}$  by a single 
$(2k+1,2k+1)$-pass-move  
where $k$ is any positive integer.

This completes the proof of Main Theorem \ref{melon}. 
\qed 
\\

\noindent
{\bf Note.}  
 \cite[Proposition 6.2 and \S6]{KauffmanNeumann},  
\cite[Proposition 5.4 and section 5]{KauffmanOgasa}, and 
\cite[section 7]{Levinecob}  imply 
Lemma \ref{Xpear} and \ref{Xbean}.

\begin{lem} \label{Xpear}  
Let $k$ be a nonnegative integer. 
Let Z=J,K. 
If $m=2k+1$, 

$
[h_{\alpha, Z\otimes^m[2]}]\cdot[h_{\beta, Z\otimes^m[2]}] 
\text{ in $V_{Z\otimes^m[2]}$}
\neq [h_{\alpha, Z}]\cdot[h_{\beta, Z}]  \text{ in $V_Z$}$
in general. 
\end{lem}

Note the difference between Lemma \ref{Xbanana} and Lemma \ref{Xpear}.

\begin{lem} \label{Xbean}
Let $k$ be a positive integer. 
$J\otimes^{2k+1}[2]$ and 
$K\otimes^{2k+1}[2]$ are non-diffeomorphic $($respectively, non-homeomorphic$)$ in general 
\end{lem}

This lemma implies  \cite[Theorem 8.3]{KauffmanOgasa} as written there: 
Let $k$ be a positive integer.  
$J\otimes^{2k+1}[2]$ and  
$K\otimes^{2k+1}[2]$ are 
 not $(2k+2, 2k+2)$-pass-move equivalent in general.

\bigbreak \noindent{\bf Note.}
In Main Theorem \ref{sugar}, 
 if we replace 
`1-links in $S^3$' 
with 
`2$\nu+$1-dimensional closed oriented submanifolds in $S^{2\nu+3}$, 
where  $\nu$ is a positive integer,',  
and 
$`(2k+1,2k+1)'$ 
with 
$`(2k+\nu+1, 2k+\nu+1)'$,  
we could prove it in the same way.  

Furthermore, in Main Theorem \ref{melon}, 
if we replace 
`2-dimensional closed oriented submanifolds in $S^4$' 
with 
`(2+2$\nu$)-dimensional closed oriented submanifolds in $S^{2\nu+4},$ 
where $\nu$ is a positive integer',  
and 
`$(2k+1,2k+2)$' 
with 
`$(2k+1+\nu, 2k+2+\nu)$',  
and 
if we remove `$k\geq2$', 
we could prove it  in the same way.  
\\

We formulate a problem.

\begin{prob}\label{Hakusan}
Let $p,q,$ and $\mu$ be positive integers. 
Suppose that 
we obtain a $(p+q-1)$-knot $J$ from a $(p+q-1)$-knot $K$ 
by one $(p,q)$-pass-move. 
Can we obtain 
$J\otimes^\mu$(Hopf) from
$K\otimes^\mu$(Hopf)  
by a $(2\mu+p, 2\mu+q)$-pass-move?
\end{prob}

\bigbreak
\section{Proof of Main Theorem \ref{coffee}}\label{bije} 

\noindent
The idea of the proof is to show that   
Farber's bijective map 
 $\mathcal S_{2m}\to\mathcal S_{2m+4} (m\geqq4)$ which is cited in 
Theorem \ref{Fsuspension} 
and a map $\mathcal S_{2m}\to\mathcal S_{2m+4} (m\geqq4)$ 
in Theorem \ref{mars}  
which is defined by using knot products 
are the same.   
In order to prove the sameness, we use the fact that each of the two maps 
is the double of a kind of suspension. 

\bigbreak

We review the classification of even dimensional simple spherical knots by Farber 
in
 \cite{Farber1978, Farber1980, Farber1983, Farber1984I, Farber1984II}.  
In particular, see \cite[\S1 and 2]{Farber1984I} for detail.

Before that, we need to define a continuous map. 
Let $n$ be a positive integer. 
Let $A$ (respectively, $B$) be a (not necessarily connected)   
$(n+1)$-dimensional compact oriented submanifold with nonempty boundary 
contained in $S^{n+2}$.  
Regard  $S^{n+2}$ as 
the one point compactification $\R^{n+2}\cup\{*\}$ of $\R^{n+2}$.  
Suppose $A\subset\R^{n+2}$, $B\subset\R^{n+2}$, and $A\cap B=\phi$.   
Take a continuous map  
$$\alpha: A\x B\to \R^{n+2}-\{0\} $$
$$\hskip-16mm (a,b)\mapsto a-b. $$

\noindent
Here we regard $a$ and $b$ as elements of $\R^{n+2}$. 
Hence we can define $a-b$ consistently. 
Thus we can define a continuous map 
$$\tau(\alpha): A\x B\to S^{n+1}$$
$$\hskip12mm (a,b)\mapsto \frac{a-b}{\vert a-b\vert}, $$

\noindent
where $\vert a-b\vert$ is the distance between $a$ and $b$. 
We can make $\tau(\alpha)$ into a continuous map 
$$\theta(\alpha): A\wedge B\to S^{n+1},$$

\noindent
where $\wedge$ denotes the smash product because 
$\tau(\alpha)\vert_{A\vee B}$ is a null-homotopic map. 
\\


Let $m$ be a positive integer. 
Let $K$ be a $2m$-dimensional simple spherical knot contained in $S^{2m+2}=\R^{2m+2}\cup\{*\}$. 
Let $V$ be an $(m-1)$-connected Seifert hypersurface for $K$. 
Push off $V$ into the positive direction of the normal bundle of $V$ in $\R^{n+2}$, 
and call it $V_+$.    
Let $A=V$ and $B=V_+$ in the definition of the above continuous map $\theta(\alpha)$.  
We abbreviate $\theta(\alpha):V\wedge V_+\to S^{2m+1}$ 
to  $\theta:V\wedge V\to S^{2m+1}$.  
By using this $\theta$, define 
an  $S$-map  
$\theta:V\wedge V \stackrel{s}{\to} S^{2m+1}$, where $\stackrel{s}{\to}$ denotes an $S$-map, 
and call it the {\it Seifert-Farber homotopy pairing} for $V$.
See also \cite{Farber1984I} for the definition of $S$-maps.  
Here we let $f:X\stackrel{s}{\to}Y$ denote an $S$-map made from a continuous map $f:X\to Y$. 
Note that we use the notation $f$ for both. 
The Seifert-Farber homotopy pairing defines a spherical pairing $(V, \theta)$. 
Let $n$ be a positive integer. 
The $n$-isometry $(V,u,z)$ corresponding to this $n$-pairing will be called the 
{\it $n$-isometry of the manifold $V$}.  
See \cite[\S1 and 2]{Farber1984I} for 
the definition of spherical pairings and that of $n$-isometries.  
Spherical $n$-pairings and $n$-isometries are in one-to-one correspondence 
by the correspondence written in \cite[\S2.4]{Farber1984I}. 
Furthermore \cite[Lemma 2.8]{Farber1984I}  shows that 
spherical $n$-pairings are contiguous if and only if the corresponding $n$-ismoetries are contiguous. See also \cite{Farber1984I} for the definition of the term `contiguous'.

\begin{thm}\label{F}{\rm {\bf(\cite[Theorem 2.6]{Farber1984I}.)}}
Let $m$ be a positive integer. 
The $R$-equivalence class of the $2m$-isometry $(V,u,z)$ 
does not depend on the choice of the Seifert hypersurface $V$ 
and is well-defined by the type of the given $2m$-dimensional simple spherical knot $K$. 
If 
$m\geqq4$, 
then the map, 
sending each $2m$-dimensional simple spherical knot 
to the $R$-equivalence class of the $2m$-isometry of some $(m-1)$-connected Seifert hypersurface of this knot, 
is a bijection of 
the set of isotopy types of $2m$-dimensional simple spherical knots 
to 
the set of $R$-equivalence classes of $2m$-isometries given on all virtual complexes of length $\leqq2$.  
\end{thm}

\smallbreak \noindent 
{\bf Note.} 
See \cite[\S1 and 2]{Farber1984I} for the definition of 
`$R$-equivalence classes'  
and 
that of 
  `$n$-isometries given on all virtual complexes of length $\leqq2$'.  

\bigbreak

We introduce a few continuous maps. 
Let $\R^{n+3}=\R\x\R^{n+2}$. 
Suppose that an element in $\R^{n+3}=\R\x\R^{n+2}$ is represented by  $(x,p)$, where $x\in\R$ and $p\in\R^{n+2}$. 
Construct the suspension $\Sigma A$    
from $A\subset\R^{n+2}$ and $(1,0), (-1,0)\in\R\x\R^{n+2}$.   
%
Of course $\Sigma A$ is not in $\{0\}\x\R^{n+2}$ but in $\R^{n+3}$. 
Define a continuous map 
$\beta_1(\alpha): \Sigma A\x B\to\R^{n+3}-\{0\}$ as follows. 
Note that $B\subset\{0\}\x\R^{n+2}\subset\R^{n+3}$. 
We use the following lemma.

\begin{lem}\label{moon}
Let $0\leqq t\leqq1$. 
Take an element in 
the intersection of  $\Sigma A\subset\R\x\R^{n+2}$ 
and 
$\{x|x=t\}\x\R^{n+2}\subset\R\x\R^{n+2}$.     
Then for an element $(0,a)\in A\subset\{x|x=0\}\x\R^{n+2}$,   
the element in the intersection 
is denoted by $(t, (1-t)a)\in\{x|x\in\R\}\x\R^{n+2}=\R\x\R^{n+2}$.
\end{lem}

Define a continuous map $\beta_{1,u}(\alpha)$ to be   
$$ \beta_{1,u}(\alpha):((\Sigma A)\cap(\R_{\ge0}\x\R^{n+2}) )\x B\to\R^{n+3}-\{0\} $$
$$\hskip39mm (p,q)\mapsto p-q. $$

\noindent
Recall $B\subset\{x|x=0\}\x\R^{n+2}$.
Hence this continuous map is represented by 

\vskip2mm\hskip50mm
$ ((t, (1-t)a), (0,b))\hskip3mm\mapsto  \hskip3mm(t, (1-t)a-b) $

\hskip93mm$=(t, a-at-b).$
\vskip2mm

\noindent in another way. 
We use the following lemma.

\begin{lem}\label{Titan}
Let $-1\leqq t\leqq0$. 
Take an element in 
the intersection of  $\Sigma A$ and \newline $\{x|x=t\}\x\R^{n+2}$.  
Then for an element $(0,a)\in A$,   
the element in the intersection 
is denoted by $(t, (1+t)a)$.  
\end{lem}

Define a continuous map $\beta_{1,d}(\alpha)$ to be  
$$ \beta_{1,d}(\alpha):((\Sigma A)\cap(\R_{\le0}\x\R^{n+2}))\x B\to\R^{n+3}-\{0\} $$
$$\hskip39mm  (p,q)\mapsto p-q. $$

\noindent 
Hence this continuous map is represented by 
$$ ((t, (1+t)a), (0,b))\mapsto  \hskip7mm(t, (1+t)a-b) $$
\hskip92mm$=(t, a+at-b).$

\noindent in another way.  
The continuous map $\beta_1(\alpha)$ is defined by using 
$\beta_{1,u}(\alpha)$ and $\beta_{1,d}(\alpha)$.

\bigbreak
We can define a continuous map 
$\beta_2(\alpha): A\x \Sigma B\to\R^{n+3}-\{0\}$ 
in the same way as we define $\beta_1(\alpha)$. 
%
Let $\beta(\alpha)=\beta_1(\beta_2(\alpha)): 
(\Sigma A)\x\Sigma B\to\R^{n+4}-\{0\}$. 
We will use \newline
$\beta(\beta(\alpha)):(\Sigma^2 A)\x\Sigma^2 B\to\R^{n+6}-\{0\}$.
Let $\beta^2(\alpha)$ denote $\beta(\beta(\alpha))$.

\bigbreak 
Define a continuous map 
$\gamma_1(\alpha):(\Sigma A)\x B\to\R^{n+3}-\{0\} $ as follows. 

\noindent 
Let $0\leqq t\leqq1$. 
Define a continuous map $\gamma_{1,u}(\alpha)$ to be   
$$\gamma_{1,u}(\alpha):((\Sigma A)\cap(\R_{\ge0}\x\R^{n+2})) \x B\to\R^{n+3}-\{0\} $$
$$ ((t, (1-t)a), (0,b))\mapsto  (t, (1-t)(a-b)).  $$

\noindent 
Let $-1\leqq t\leqq0$. 
Define a continuous map $\gamma_{1,d}(\alpha)$ to be   
$$ \gamma_{1,d}(\alpha):((\Sigma A)\cap(\R_{\le0}\x\R^{n+2}))\x B\to\R^{n+3}-\{0\} $$
$$ ((t, (1+t)a), (0,b))\mapsto  (t, (1+t)(a-b)).  $$

\noindent 
The continuous map $\gamma_1(\alpha)$ is defined by using 
$\gamma_{1,u}(\alpha)$ and $\gamma_{1,d}(\alpha)$.

\bigbreak
We can define a continuous map 
$\gamma_2(\alpha): A\x \Sigma B\to\R^{n+3}-\{0\}$ 
 in the same way as we make $\gamma_1(\alpha)$. 
%
Let $\gamma(\alpha)=\gamma_1(\gamma_2(\alpha)): 
(\Sigma A)\x\Sigma B\to\R^{n+4}-\{0\}$. 
We will use \newline
$\gamma(\gamma(\alpha)):(\Sigma^2 A)\x\Sigma^2 B\to\R^{n+6}-\{0\}$.
Let $\gamma^2(\alpha)$ denote $\gamma(\gamma(\alpha))$.

\begin{lem}\label{sun}
The above continuous maps $\beta_1(\alpha)$ and $\gamma_1(\alpha)$ are homotopic. 
\end{lem}

\noindent 
{\bf Proof of Lemma \ref{sun}.}   
We define a homotopy  
$ \xi_s:\Sigma A \x B\to\R^{n+3}-\{0\},$ where $0\leq s\leq1$,   
such that $\xi_0=\beta$ and that $\xi_1=\gamma$ as follows.

\noindent 
Let $0\leqq t\leqq1$. 
Define a continuous map $\xi_{s,u}$ to be   
$$ \xi_{s,u}:(\Sigma A)\x B\to\R^{n+3}-\{0\} $$
$$ \hskip40mm((t, (1-t)a), (0,b)) \hskip3mm \mapsto \hskip3mm
(t, s(1-t)(a-b)+(1-s)(a-b-ta))$$
\hskip88mm$=(t, (1-t)a-(1-st)b).$

\noindent 
Let $-1\leqq t\leqq0$. 
Define a continuous map $\xi_{s,d}$ to be   
$$\xi_{s,d}:(\Sigma A)\x B\to\R^{n+3}-\{0\} $$

\vskip2mm\hskip40mm
$((t, (1+t)a), (0,b))\hskip3mm\mapsto  \hskip3mm(t, s(1+t)(a-b)+(1-s)(a-b+ta))$

\hskip88mm$=(t, (1+t)a-(1+st)b).$
\vskip2mm

The continuous map $\xi_s$ is defined by using $\xi_{s,u}$ and $\xi_{s,d}$.  
We claim that $\xi_{s,u}$, $\xi_{s,d}$, and $\xi_s$ are well-defined.   
{\it Reason}:
Let $p\in\Sigma A \x B$. 
Let $\xi(p)=(t,\kappa)$. 
If $t\neq0$, $\xi_s(p)\in\R^{n+3}-\{0\}.$  
If $t=0$, then $\kappa=a-b$. 
Hence $\kappa\neq0$. Hence $\xi_s(p)\in\R^{n+3}-\{0\}.$  \qed


\bigbreak 
We can prove the following lemma in the same way as we prove Lemma \ref{sun} 
\begin{lem}\label{summer}
The above continuous maps $\beta^2(\alpha)$ and $\gamma^2(\alpha)$ 
are homotopic. 
\end{lem}

\bigbreak 
We state a theorem which was proved in  \cite{Farber1978, Farber1980, Farber1983, Farber1984I, Farber1984II} after some preparations 
(see Theorem \ref{Fsuspension}).  
Take the suspension $\Sigma X$ of $X$. 
$\Sigma X$ is made from $X$ and $[-1,1]$ as usual. 
An element of $\Sigma X$ is represented by $(t,x)$ by using  $t\in[-1,1]$ and $x\in  X$ as usual. This way of representation is different from the previous way that we used in order to define $\beta(\alpha)$ and $\gamma(\alpha)$.

Take a continuous map 
$$f:X\to Y$$
$$\hskip5mm x\mapsto y. $$

\noindent
 Define a continuous map to be  
$$\Sigma f:\Sigma X\to\Sigma Y$$
$$\hskip8mm (t,x)\mapsto(t,y),$$
and call it the {\it suspension} of $f$. 
Note that 
$(1,x)$ 
(respectively, $(-1,x)$)   
and 
$(1,x')$ (respectively, $(-1,x')$)  
represent a same element even if $x\neq x'$.

\noindent
Take a continuous map 
$$f:A\x B\to Y$$
$$\hskip5mm (a,b)\mapsto y. $$
Define a continuous map to be  
$$\Sigma^{1,0} f:(\Sigma A)\x B\to\Sigma Y$$
$$\hskip15mm ((t,a), b)\mapsto(t,y),$$
and call it the {\it $(1,0)$-suspension} of $f$. 
Define a continuous map to be  
$$\Sigma^{0,1} f:A\x\Sigma B\to\Sigma Y$$
$$\hskip15mm (a,(t,b))\mapsto(t,y),$$
and call it the {\it $(0,1)$-suspension} of $f$. 
We can define $\Sigma^{\nu,\mu} f$ in the same manner. 
We have the following.

\begin{lem}\label{key}
Let $n$ be a positive integer. 
Take $\alpha$, $\tau(\alpha)$ as above. 
Then $\Sigma^{2,2}\tau(\alpha):\Sigma^2 A\x \Sigma^2 B\to \Sigma^4 S^{n+1}$  
is homotopic to 
$\tau(\gamma^2(\alpha))$.  
$($Note: $\Sigma^4 S^{n+1}$ is regarded as $S^{n+5}.)$  
\end{lem}  

Furthermore we have the following.

\begin{lem}\label{water}  
Let $\psi,\phi,\nu, \nu-\psi,$ and $ \nu-\phi$ be nonnegative integers. 
$(\Sigma^\psi A)\wedge  \Sigma^{\nu-\psi}B$ is homotopy type equivalent to  
$(\Sigma^\phi A)\wedge  \Sigma^{\nu-\phi}B$.  
\end{lem}   

Make
$$\Sigma^{\psi, \nu-\psi}f:(\Sigma^\psi A)\x  \Sigma^{\nu-\psi}B\to Y$$
into a continuous map  
$$(\Sigma^\psi A)\wedge  \Sigma^{\nu-\psi}B\to Y.$$ 
Make
$$\Sigma^{\psi, \nu-\phi}f:(\Sigma^\phi A)\x  \Sigma^{\nu-\phi}B\to Y$$
into a continuous map 
$$(\Sigma^\phi A)\wedge  \Sigma^{\nu-\phi}B\to Y.$$ 
By Lemma \ref{water} these two new continuous maps are regarded as the same one,  and call it  $\Sigma^{\nu} f.$  
Hence Lemma \ref{key} implies the following.

\begin{lem}\label{important}
$\Sigma^4\theta(\alpha):(\Sigma^2 A)\wedge \Sigma^2 B\to S^{n+5}$  
is homotopic to  \newline
$\theta(\gamma^2(\alpha)):(\Sigma^2 A)\wedge \Sigma^2 B\to S^{n+5}$.  
\end{lem}

\noindent
{\bf Note.}  
Recall that for finite CW complexes, 
the reduced suspension of X is homotopy type equivalent to the ordinary suspension.
The suspensions are used in \cite[\S1 and 2]{Farber1984I} and is also used in 
the following Theorem \ref{Fsuspension}.


\begin{thm}\label{Fsuspension}   {\rm{\bf(\cite[Theorem 2.6 and \S1 and 2]{Farber1984I}.)}}
Under the condition $(*)$, we have the conditions $(1)$ and $(2)$. 
 
\smallbreak \noindent$(*)$  
Let $m$ be an integer$\geqq4$.    
Let $K^{2m}$ be a simple spherical $2m$-knot. 
Let $V^{2m+1}$ be an $(m-1)$-connected Seifert hypersurface for $K^{2m}$. 
Let $\theta: V^{2m+1}\wedge V^{2m+1}\stackrel{s}{\to}S^{2m+1}$ 
be a Seifert-Farber homotopy pairing. 

\smallbreak \noindent$(1)$  
Make an $S$-map  
$\Sigma^4\theta:(\Sigma^2V^{2m+1})\wedge\Sigma^2V^{2m+1}\stackrel{s}{\to}S^{2m+5}$ 
by using the suspensions.  
Then there is  a simple spherical $(2m+4)$-knot $J^{2m+4}$  $\subset S^{2m+6}$   
with the following property:     
Let $U^{2m+5}$ be an $(m+1)$-connected Seifert hypersurface for $J^{2m+4}$. 
Let 
$\rho:U^{2m+5}\wedge U^{2m+5}\newline\stackrel{s}{\to} S^{2m+5}$ 
be a Seifert-Farber homotopy pairing. 
Then $\Sigma^4\theta$ is $R$-equivalent to $\rho$. 

\smallbreak \noindent$(2)$    
For each integer $m\geqq4$, 
the operation sending $K^{2m}$ to  $J^{2m+4}$ in $(1)$ gives a one-to-one map 
$\mathcal S_{2m}\to\mathcal S_{2m+4}$.   
\end{thm}

Lemma \ref{important} implies the following. 

\begin{lem}\label{gohan} 
Suppose the condition $(*)$ in Theorem $\ref{Fsuspension}$. 
Suppose that $\theta$ is made from a continuous map 
$$\alpha_V: V^{2m+1}\x V^{2m+1}_+\to \R^{2m+2}-\{0\}.$$
$$\hskip-1mm (a,b)\mapsto a-b. $$
That is, $\theta$ is an $S$-map made from $\theta(\alpha_V)$. 
$($This $S$-map is also called $\theta(\alpha_V).)$  

Then we have the following: $\theta(\gamma^2(\alpha_V))$ is $R$-equivalent to $\rho$ 
in Theorem $\ref{Fsuspension}.(1)$,  
and is a Seifert-Farber homotopy pairing for $J^{2m+4}$. 
The correspondence $\theta(\alpha_V)\to\theta(\gamma^2(\alpha_V))$ gives the bijection 
$\mathcal S_{2m}\to\mathcal S_{2m+4}$  
for each integer $m\geqq4$.    
\end{lem}


\begin{thm}\label{mars} {\rm{\bf (\cite[Lemma 6.1, its proof and section 6]{KauffmanNeumann}.)}}
Under the condition $(\natural)$, we have the condition $(\#)$.   

\smallbreak\noindent  $(\natural)$ 
Let $m$ be an integer$\geqq2$. 
Let $K^{2m}$ be a simple spherical $2m$-knot in $S^{2m+2}\\=\R^{2m+2}\cup\{*\}$.   
Suppose that $K^{2m}\subset\R^{2m+2}.$ 
Let $V^{2m+1}$ 
be an $(m-1)$-connected Seifert hypersurface contained in $\R^{2m+2}$ for $K^{2m}$. 
Make $V^{2m+1}_+$ from $V^{2m+1}$ as usual. 
Take a continuous map 
$$\alpha_V: V^{2m+1}\x V^{2m+1}_+\to \R^{2m+2}-\{0\}$$
$$\hskip-1mm (a,b)\mapsto a-b. $$
Take $\beta^2(\alpha_V)$ as above. 
Make $K^{2m}\otimes{\rm Hopf}$ in $S^{2m+6}=\R^{2m+6}\cup\{*\}$.

\smallbreak\noindent  $(\#)$ 
There is an $(m+1)$-connected Seifert hypersurface $W^{2m+5}$ for 
$K^{2m}\otimes{\rm Hopf}$ with the following properties.

\smallbreak\noindent  $(1)$ 
$W^{2m+5}$ is homotopy type equivalent to $\Sigma^2 V^{2m+1}$ 
by continuous maps \newline
$\zeta: \Sigma^2 V^{2m+1}\to W^{2m+5}$ and 
$\zeta': W^{2m+5}\to\Sigma^2 V^{2m+1}$.

\noindent
Make $W^{2m+5}_+$ from $W^{2m+5}$ as usual.    
Then 
$W^{2m+5}_+$ is homotopy type equivalent to $\Sigma^2 V^{2m}_+$ 
by continuous maps 
$\zeta_+: \Sigma^2 V_+^{2m+1}\to W_+^{2m+5}$ 
and $\zeta_+': W_+^{2m+5}\to\Sigma^2 V_+^{2m+1}$.

\smallbreak
\noindent$(2)$  
Take a continuous map 
$$\alpha_W: W^{2m+5}\x W^{2m+5}_+\to \R^{2m+6}-\{0\}$$
$$\hskip-1mm (a,b)\mapsto a-b. $$
Take the continuous map 
$\zeta\x\zeta_+:(\Sigma^2V^{2m+1})\x\Sigma^2V^{2m+1}_+\to W^{2m+5}\x W^{2m+5}_+$. 
Then the continuous map $\alpha_W\circ(\zeta\x\zeta_+)$ is homotopic to 
the continuous map $\beta^2(\alpha_V)$. 
\end{thm}

The $S$-map made from $\theta(\alpha_W\circ(\zeta\x\zeta_+))$ 
is a Seifert-Farber homotopy pairing for $K^{2m}\otimes{\rm Hopf}$. 
By Theorem \ref{F} we have the following.

\begin{lem}\label{awa} 
Under the condition $(\natural)$ in Theorem $\ref{mars}$,  
 the $S$-map made from $\theta(\beta^2(\alpha_V))$ 
is a Seifert-Farber homotopy pairing for $K^{2m}\otimes{\rm Hopf}$. 
\end{lem}

Lemmas \ref{summer} and \ref{awa} implies the following.

\begin{lem}\label{mugi} 
Under the condition $(\natural)$ in Theorem $\ref{mars}$,  
 the $S$-map made from $\theta(\gamma^2(\alpha_V))$ 
is a Seifert-Farber homotopy pairing for $K^{2m}\otimes{\rm Hopf}$. 
\end{lem}

Lemmas \ref{gohan} and \ref{mugi} imply Main Theorem \ref{coffee}.   
This completes the proof of Main Theorem \ref{coffee}.\qed

\bigbreak \noindent{\bf Note.}  
We could extend Main Theorem \ref{coffee} to the case of the stable knots corresponding to 
the set of $R$-equivalence classes of $2m$-isometries given on all virtual complexes of 
a given length.

\bigbreak\section{Problem}\label{problem}
Let $k$ be a positive integer. 
One way of saying, 
investigating 
how far $\otimes^k$Hopf$:C_1\\\to C_{4k+1}$ and $\otimes^k$Hopf$:\mathcal K_1\to \mathcal S_{4k+1}$ are different from  isomorphism maps    
has generated many fruitful results on 1-knots 
as the readers know.
It is natural to formulate the following problem. 

\begin{prob}\label{tsugi}
Let $\mu$ be an integer$\geqq2$. 
Investigate $\otimes^\mu$Hopf$:\mathcal K_2\to \mathcal S_{4\mu+2}$.   

Let $\nu$ be an integer$\geqq1$.  
Investigate $\otimes^\nu$Hopf$:\mathcal S_4\to \mathcal S_{4\nu+4}$.   

Let $\tau$ be an integer$\geqq1$.  
Investigate $\otimes^\tau$Hopf$:\mathcal S_6\to \mathcal S_{4\tau+6}$.   

Investigate $\otimes$Hopf$:\mathcal K_2\to \mathcal S_{6}$.  
\end{prob}

It will be exciting to combine Problems \ref{asu} and \ref{tsugi}. 
Make problems from Problem \ref{asu} replacing the ribbon-move by other local moves. 
Combine the problems with  Problem \ref{tsugi}. It will be also exciting.


\bigbreak \noindent
Louis H. Kauffman: 
Department of Mathematics, Statistics, and Computer Science,  
University of Illinois at Chicago, 
851 South Morgan Street, 
Chicago, Illinois 60607-7045, USA  \quad
kauffman@uic.edu

\bigbreak \noindent
Eiji Ogasa:  Computer Science, Meijigakuin University, Yokohama, Kanagawa, 244-8539, Japan 
\quad pqr100pqr100@yahoo.co.jp  \quad
ogasa@mail1.meijigkakuin.ac.jp

\end{document}